\documentclass[11pt]{article}
\usepackage{amssymb}
\usepackage{graphicx}
\usepackage{amsmath}
\usepackage{multirow}
\usepackage{caption}
\usepackage{subcaption}
\usepackage{enumerate}
\usepackage{mathrsfs}
\usepackage[labelfont=bf]{caption}

\newcommand{\R}{\mathbb{R}}

\newcommand{\om}{\omega}
\DeclareMathOperator{\I}{I}
\DeclareMathOperator{\II}{II}
\DeclareMathOperator{\III}{\cal{PD}}
\DeclareMathOperator{\Hess}{Hess}

\DeclareMathOperator{\Ind}{Ind}

\newtheorem{theorem}{Theorem}[section]
\newtheorem{proposition}[theorem]{Proposition}
\newtheorem{lemma}[theorem]{Lemma}
\newtheorem{remark}[theorem]{Remark}
\newtheorem{corollary}[theorem]{Corollary}
\newtheorem{definition}[theorem]{Definition}
\newtheorem{examples}[theorem]{Example}

\begin{document}
\title{Umbilic Points on the Finite and Infinite Parts of Certain Algebraic Surfaces}

\author{Brendan Guilfoyle\thanks{School of STEM, Munster Technological University, Kerry,
		Tralee, Co. Kerry, Ireland \newline e-mail: brendan.guilfoyle@ittralee.ie} \ 
	and Adriana Ortiz-Rodr\'iguez\thanks{Instituto de Matem\'aticas, Universidad 
		Nacional Aut\'onoma de M\'exico. Area de la Inv. Cient., Circuito exterior C.U.,
		 Mexico City 04510, M\'exico \newline e-mail: aortiz@matem.unam.mx}}
\date{}
\maketitle

\begin{abstract}
\noindent The global qualitative behaviour of fields of principal directions for the 
graph of a real valued polynomial function $f$ on the plane is studied. 
We determine and analyze the projective extension of these fields and show that it 
is defined by an analytic quadratic form on the whole unit 2-sphere. We prove 
that every umbilic point at infinity of this extension has a Poincar\'e-Hopf index 
equal to 1/2 and the topological type of a Lemon or a Monstar. 
As a consequence of these results we provide a Poincar\'e-Hopf type formula for the graph of $f$ 
pointing out that, if all umbilics are isolated, the sum of all indices of the principal directions 
at its umbilic points only depends upon 
the number of real linear factors of the homogeneous part of highest degree of $f$. 
A similar analysis is carried out in the case that $f$ is a homogeneous polynomial.
\end{abstract}\medskip



\vspace{0.1in}

\section{Introduction}

\begin{figure}[h!]    
	\begin{center}
		\includegraphics[width=3.60in]{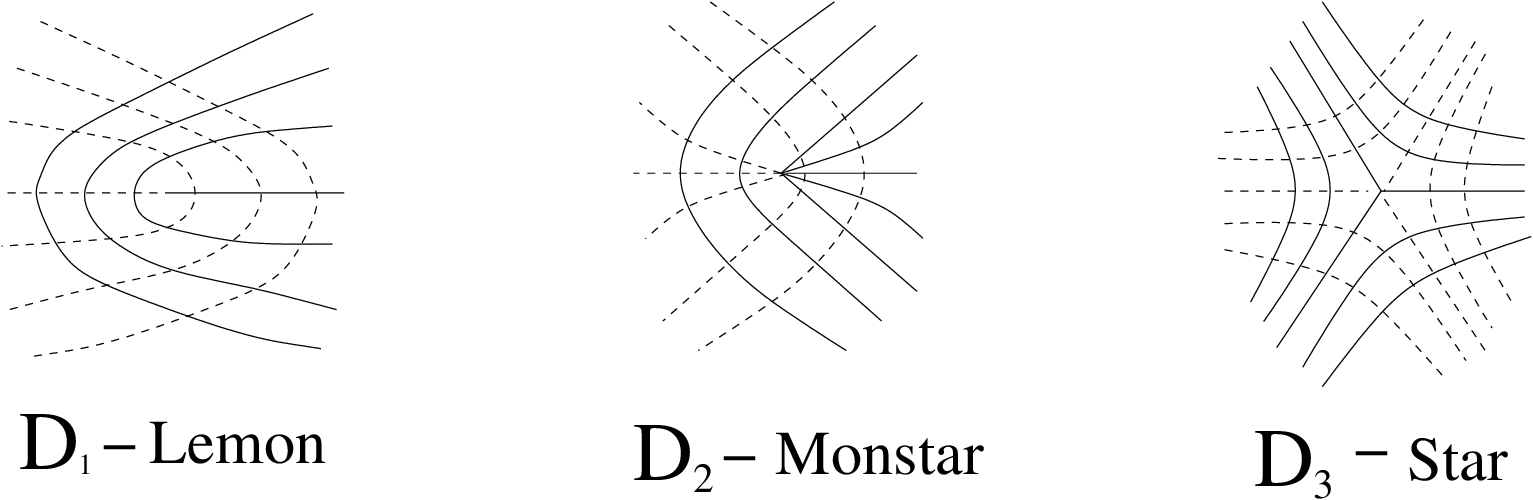}
		\caption{Darbouxian points}
		\label{darboux-points} 
	\end{center}
\end{figure}

For any oriented smooth surface in real Euclidean 3-space, the eigenspaces of its second 
fundamental form define two orthogonal line fields, called {\it fields of principal directions}, 
whose singularities are the {\it umbilics of the surface}. The study of the fields of principal 
directions and the principal lines of a smooth 
surface dates back to Euler, Darboux \cite{darb}, Monge \cite{monge} and Cayley \cite{Cay}, 
amongst others. An umbilic is 
characterized by the fact that its principal curvatures are equal. Moreover, to each isolated 
umbilic can be attached the index of either one of the two fields. This index is of the form 
${n}/{2}$, with $n \in \mathbb{Z}$. When such a surface is generic, the behaviour of the 
principal lines in the neighbourhood of an umbilic can only be one of three {\it Darbouxian types}: 
Lemon, Monstar and Star \cite{soto-gut} (see Fig. \ref{darboux-points}) with Poincar\'e-Hopf index 
$\frac{1}{2}, \frac{1}{2}, -\frac{1}{2}$, respectively. For topological reasons, if all of the 
umbilics are isolated, the sum over all half-integer indices of the fields of principal directions 
at its umbilic points equals the Euler characteristic of the surface.
\medskip

The study of umbilic points at infinity of a surface has been developed from various perspectives previously.
In the case of smooth surfaces, V. Toponogov analyzes in \cite{Top}, 
surfaces $S$ homeomorphic to a plane which are complete and convex. In terms of the principal curvatures 
 $\kappa_1,\kappa_2$ of $S$ he states the conjecture:

{\centerline{\it ``on a complete convex surface $S$ homeomorphic to a plane, the 
equality}}
{\centerline{\it inf$_{p \in S}$ $| k_2(p) - k_1(p)| = 0$ holds".}}\medskip

\noindent If there are no finite umbilic points, this can be taken to mean that 
there must be an umbilic point {\em at infinity}. In the same paper, he proves 
this conjecture under some additional hypotheses. We remark that after a 
straightforward calculation, this equality can be verified for any surface that is 
given as the graph of a real polynomial.\medskip

Another instance of umbilic points at infinity in the smooth case is the study, 
carried out by R. Garcia and J. Sotomayor in \cite{Gar-Sot-06}, on stable patterns 
of the nets of principal curvature lines on surfaces embedded in Euclidean 3-space 
near their end points, at which the surfaces tend to infinty. The research just cited is 
an extension of the work by the same authors \cite{Gar-Sot-96} 
and devoted to the analysis {\it at infinity} of the principal curvature 
nets of smooth algebraic surfaces in real Euclidean 3-space.  
The surfaces discussed in \cite{Gar-Sot-96} do not cover those studied in this 
paper as the former consider surfaces having a smooth projective closure while 
those studied in this paper are singular at infinity.\medskip

In the particular case of a surface given by the graph of a homogeneous polynomial 
$f\in \mathbb{R}[x,y]$, 
the study of the index of the umbilic point that appears after the one-point 
compactification of the surface  with the point at infinity, was carried out
by N. Ando in \cite{ando}.  \medskip

On the other hand, and in a broader context, there is the investigation of 
singular points, their Poincar\'e-Hopf index and topological type, that arise {\it at 
infinity} as a result of the projective extension of a quadratic differential form 
on the plane. In {\rm \cite{Gui}} V. Gu\'{\i}\~nez considers the set ${\cal F}_m$ 
of positive quadratic forms 
\begin{equation*}
\om = a(x,y)\, dy^2 + b(x,y)\, dx dy + c(x,y)\, dx^2,
\end{equation*}
such that $a, b, c \in \R[x,y]$ are polynomials of degree at most $m$, the 
function  $\, b^{2}-4ac\,$ is non-negative at every point of the $xy$-plane, and 
$(b^2-4ac)^{-1}(0) = a^{-1}(0) \cap b^{-1}(0) \cap c^{-1}(0)$. 
For a generic form $\om \in {\cal F}_m$, he studies the projective extension of 
$\om$ and proves, amongst other things, that the topological behavior of these 
foliations in a neighborhood of a {\it singular point at infinity} is a Monstar or 
a Star (Remark 2.9 of \cite{Gui}).	\medskip

In this paper, we study the global qualitative behaviour of fields of principal 
directions for the graph of a polynomial $f \in \R[x,y]$. This study is carried 
out through the analysis of the projective extension of the quadratic form $\III$
that defines the fields of principal directions of the surface,  
and the singular points that appear on it. It is worth 
emphasizing that even though the quadratic form $\III$ belongs to the set 
${\cal F}_{3n-4}$, it is not a generic form of those analyzed in {\rm \cite{Gui}} 
as the polynomial expression determining the singular points at infinity of the 
projective extension of $\III$ vanishes at every point on the equator of the 
unit sphere. Nevertheless, by using Euler's Lemma we obtain a  
non-degenerate form on the equator. \medskip

In what follows we begin by providing an analytic quadratic 
form $\Phi$ defined on the whole sphere that describes the projective extension of 
the form $\III$. The two solution fields of this form, $ \mathbb{Y}_{1},\mathbb{Y}_{2}$, 
being restricted to the upper or lower open hemispheres, are diffeomorphic to the fields 
of principal directions, in Theorem \ref{extendedqde}. 
We then study the topological properties of the isolated singular points on the 
equator of these solution fields. In Theorem \ref{indice-umbilic-inf}, we prove that 
every such singular point, called {\it an umbilic point at infinity} has a Poincar\'e-Hopf 
index equal to $1/2$ and the topological type of a Lemon if $n=2$ and for $n\geq 3$, of a Monstar. 
\medskip

Regarding the number of umbilic points at infinity that can appear in the projective 
extension of the form $\III$, in Theorem \ref{number-umbilic-inf} we establish an upper 
bound. In Corollary \ref{ellip-hyperb-cases} we analyze the following 
two special cases.
A homogeneous polynomial on $\mathbb{R}[x,y]$ is {\it elliptic} ({\it hyperbolic})
if its Hessian function has no real linear factors and it is non-negative (non-positive). 
When  $f_n$, the highest degree homogeneous part of a polynomial $f$ of degree $n$, is 
elliptic it is proven that there are no umbilic points at infinity. 
If $f_n$ is 
hyperbolic, the number of umbilic points at infinity is bounded by twice the number of 
real linear factors of $f_n$.\medskip

In subsection \ref{remarks-caso-homog} we prove some remarks of the homogeneous case.
In  Theorem \ref{umbilic-inf-hom} we prove that when $f$ is a homogeneous polynomial 
any flat point on the equator is an umbilic point at infinity. 
\medskip

One of the main results of this paper is Theorem \ref{poincare-formula}, which provides 
a Poincar\'e-Hopf type formula for the graph of $f$. This shows that the sum of the 
indices over all umbilic points only depends upon the number of real linear factors of $f_n$. \medskip

In section \ref{examples} we display the global configurations of the fields $\mathbb{Y}_1, 
\, \mathbb{Y}_2$ for some specific cases.
We conclude the paper with the proof of Theorem \ref{indice-umbilic-inf} developed 
in section \ref{sect-6}.

\section{Preliminaries} 

This section provides some definitions and basic results that will be essential 
in the rest of the article. 
In section \ref{subsection-2.1}, we define the differential form of principal directions 
$\,\III$ of the graph of a differentiable function $f: \R^2 \rightarrow \R$ which will 
be used to determine the projective extension of the fields of principal directions of $f$.
In section \ref{determ-campos-globalm}, we recall the fact that any smooth positive 
quadratic differential form defined on an orientable smooth surface determines globally 
two direction fields. In section \ref{subsection-2.3}, we define the projective Hessian 
curve of a polynomial $f \in \R[x,y]$ which will be used in Theorem \ref{indice-umbilic-inf}.

\subsection{Fields of Principal Directions}\label{subsection-2.1}
Given a smooth surface $S$ in Euclidean 3-space the Gauss 
map $\, N : S \rightarrow \mathbb{S}^2$ associates a unit vector normal to a point $p$ on $S$ 
in a smooth way, as long as $S$ is orientable. The eigenvalues $\, -k_1, -k_2$ 
of the operator $\, D N|_p : T_pS \rightarrow T_pS$ define the {\it principal 
curvatures} $\, k_1, k_2$ {\it of the surface at the point $p$}. The points on 
$S$  at which the principal curvatures coincide are called {\it umbilic points}. 
For any non-umbilic point $p$ on $S$ the 
eigenspaces of $DN|_p$, associated to $\, -k_1$ and $-k_2$ are two orthogonal 
directions on $T_pS$ called {\it principal directions}. These directions determine two smooth direction fields which are mutually orthogonal. The maximal 
integral curves of the fields of principal directions are 
called {\it the principal curvature lines of the surface}.\medskip

In order to understand the global behaviour of the principal curvature 
lines on the graph of a differentiable function $f: \R^2 \rightarrow \R$, it is 
useful to consider the projection map
$\,\pi : \mathbb{R}^3 \rightarrow \mathbb{R}^2 $, $(x,y,z)\mapsto (x,y)$.
The image under $\pi$ of the fields of principal directions yields two fields of 
lines that are described by the quadratic differential equation:
\begin{equation}\label{gral-principal-fields}
(Eq-eQ) dx^2 + (Eg-eG) dx dy + (Qg-qG) dy^2 = 0 ,
\end{equation}
where
$\hskip 2.4cm \I (x,y) = E(x,y) dx^2 +2 Q(x,y) dxdy + G(x,y)dy^2,$ \medskip 

$\hskip 2.9cm \II (x,y) = e(x,y) dx^2 +2 q(x,y) dxdy + g(x,y) dy^2$\medskip

\noindent are, respectively the first and second fundamental forms of the surface. 
A point on the $xy$-plane is the projection of an umbilic point on $S$ if and only 
if the coefficients of the form in equation (\ref{gral-principal-fields}) vanish 
at such point. \medskip

After a direct simplification of the coefficients of equation (\ref{gral-principal-fields}), it becomes
\begin{align}\label{principal-fields-eqtn}
\Big( f_{xy} +f_{xy}(f_x)^2 &- f_x f_y f_{xx}\Big) dx^2 +
\Big( f_{yy} \Big(1+(f_x)^2\Big)-  f_{xx} \Big(1+(f_y)^2\Big)\Big) dx dy  \nonumber\\
& + \Big( f_x f_y f_{yy} - f_{xy} - f_{xy} (f_y)^2 \Big) dy^2 = 0.
\end{align}
The differential form on the left side of equation
(\ref{principal-fields-eqtn}) is called {\it the form of principal directions} and 
will be denoted by $\,\III$. The two fields at which it vanishes will be denoted 
$\mathbb{X}_{1}$ and $\mathbb{X}_{2}$.
For the sake of simplicity, we identify the principal direction fields on the graph 
of $f$ with the fields $\mathbb{X}_{1}$ and $\mathbb{X}_{2}$.\medskip

In the particular case of an $n$-degree polynomial $f\in \mathbb{R}[x,y]$, the 
coefficients of the form $\III$ are also polynomials in $\R[x,y]$ of degree at most $3n-4$. 


\subsection{Global Determination of a Direction Field}\label{determ-campos-globalm}
\begin{figure}[h!]    
	\hskip 6cm	\includegraphics[width=3.40cm]{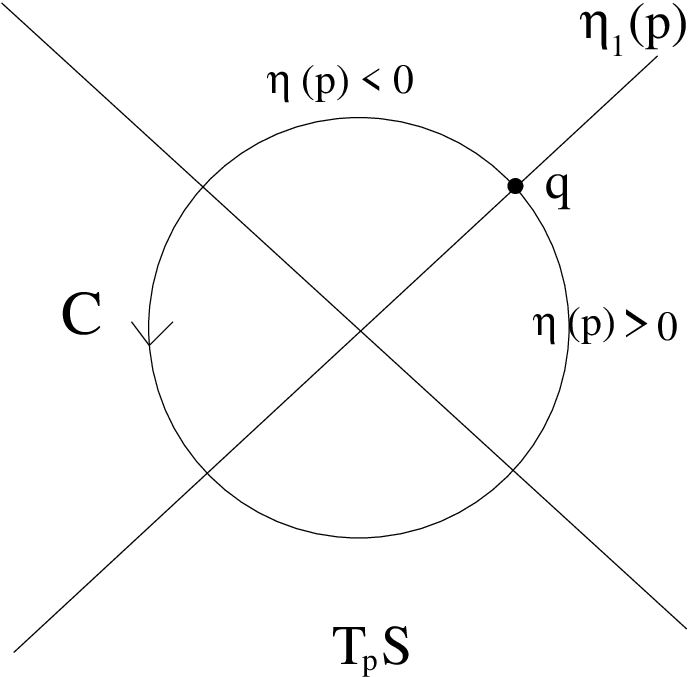}
		\caption{Determination of a direction field}
		\label{determinacion-campo-direc} 
\end{figure}

Let $S$ be an oriented smooth surface embedded in Euclidean space and $\eta$ be 
a quadratic differential form defined on an open subset $U \subset S$. Let $\eta (p): T_pS 
\rightarrow \R$ the quadratic form obtained by the restriction of $\eta$ to $T_pS$.
We say that 
{\it $\eta$ is positive} if for every point $p \in U$ the subset $\eta 
(p)^{-1}(0)$ of $T_pS$ is either the union of two transversal directions or all 
$T_pS$. 
When $\eta (p)^{-1}(0) = T_pS$ we say that {\it $p$ is a singular point of 
$\eta$}. Assume that  $p$ is not a singular point of $\eta$ and consider one 
of the two directions in $\eta (p)^{-1}(0)$ \cite{Gut-Gui}. 
Choose an oriented circle $C$ on $T_p S$ 
whose center is the origin and denote by $q$ an intersection point of $C$ with the 
chosen direction (see Fig. \ref{determinacion-campo-direc}). Consider an oriented small 
arc $\mathscr{C} = (q_1, q_2)$ on $C$ (according to the orientation of $C$) that 
contains the point $q$. Denote the chosen direction by $\eta_1 (p)$ if the 
form $\eta (p)$ is positive along the subarc $(q_1, q)$ and, negative on the subarc 
$(q, q_2)$. Otherwise, denote the chosen direction as $\eta_2 (p)$.
In this way, the set of directions $\eta_1 (p) $ obtained by varying $p$ in $U$, 
determines a continous direction field tangent to $S$.\\

\subsection{Homogenization of the Hessian Function of $f$}\label{subsection-2.3}

The projection of the parabolic curve on the graph of a smooth function $f: U\subset \R^2 \rightarrow \R$ under 
the map $\pi$ is the zero locus of the determinant of the Hessian matrix of $f$.
This determinant, $|\Hess f|=f_{xx}f_{yy}-f_{xy}^{2}$, will be referred to as the
{\it Hessian function of $f$} and its zero locus will be called the {\it Hessian 
curve of $f$}. The hyperbolic and elliptic domains are projected under $\pi$ 
into sets on which the Hessian function of $f$ is negative and positive, respectively. \medskip

When $f$ is a polynomial of degree $n$ in $\mathbb{R}[x,y]$, its 
Hessian curve is thus a real plane algebraic curve of degree at most $2n-4$. 
Considering {\it the homogeneous decomposition of $f$}, 
$f =\sum_{i=0}^{n}f_{i}$, where $f_i\in \mathbb{R}\left[ x,y\right]$ 
is a homogeneous polynomial of degree $i$, we remark that 
$$ |\Hess f| = \displaystyle\sum_{j=0}^{2n-4} h_{j},
\mbox{ where } h_j \mbox{ is a homogeneous polynomial of degree $j$ and } h_{2n-4} = |\Hess f_n|.$$
\begin{definition}\label{d:projhess}
The projective Hessian curve of $f$ {\rm is the zero locus of the 
homogeneous polynomial $H_f \in \mathbb{R}[x,y,z]$, the homogenization 
of the polynomial $|\Hess f(x,y)|$.}
\end{definition}
It follows, from the homogeneous decomposition of $f$, that $H_f$ has the expression:
$\,H_f \left(x,y,z\right) =
\sum_{j=0}^{2n-4} z^{2n-4}\, h_{j}\left(\frac{x}{z},\frac{y}{z}\right)$. 
Thus, the restriction of $H_f$ to the line at infinity $z=0$ is 
$$H_f(x,y,0) =  |\Hess f_n(x,y)|.$$ 

 
\section{Projective Extension}\label{non-hom-case}
In this section, $f $ will denote a polynomial in $\mathbb{R}[x,y]$. So, the coefficients  
of the form of principal directions $\III$ defined in the left-side of equation 
(\ref{principal-fields-eqtn}) are polynomials of degree at most $3n-4$.
Our goal is to study the fields of principal directions at infinity and we will do this 
through the projective extension of the quadratic form $\III$. 
We start by developing the extension to infinity of the form of principal directions $\III$. 
This extension will be called {\it the projective extension of $\III$} and will be obtained through 
the so-called {\it projection into the Poincar\'e sphere} \cite{poinc}. \medskip

\begin{figure}[h!] 
	\begin{center}
		\qquad\quad\includegraphics[width=4.83cm]{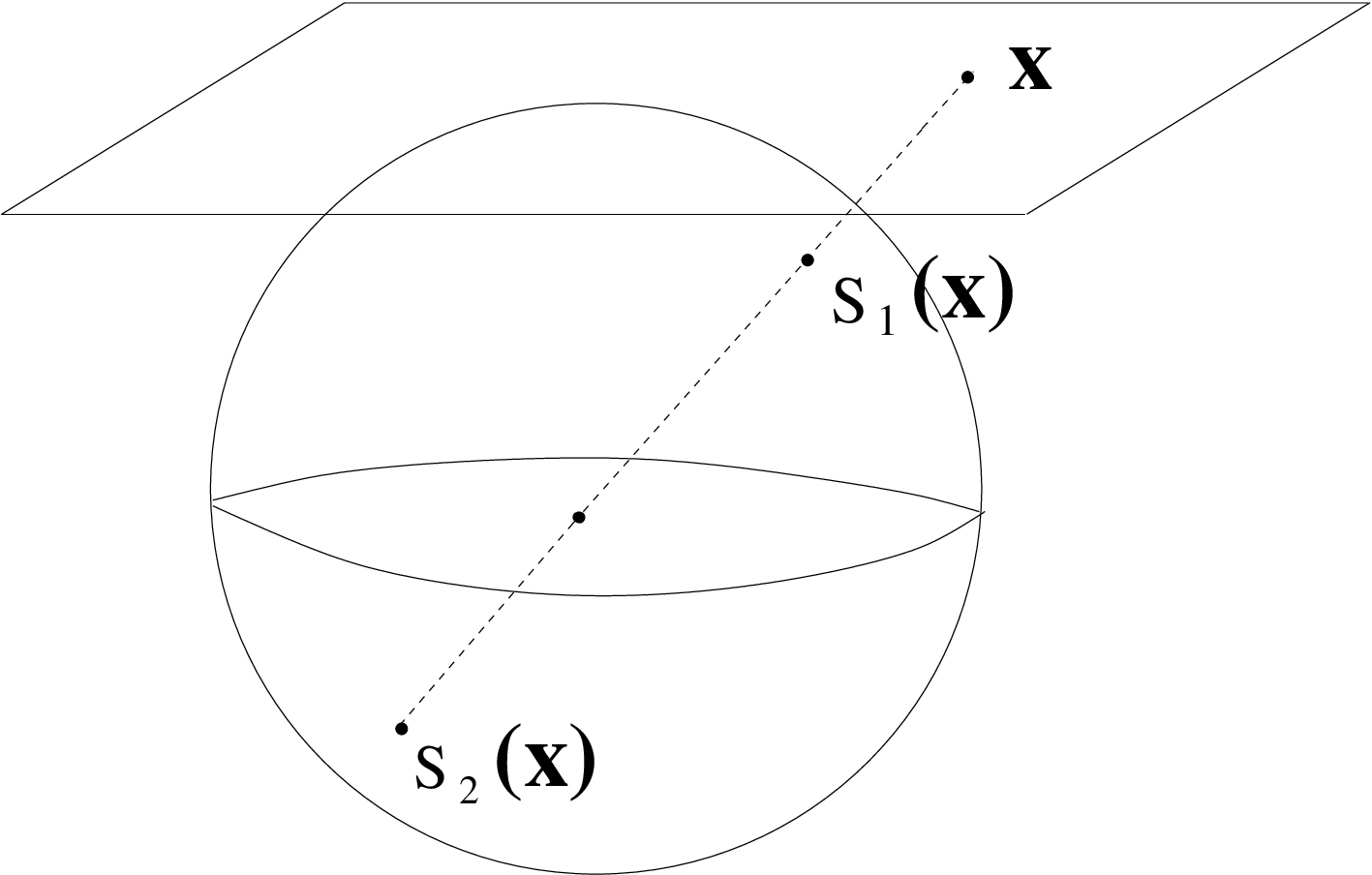}
		\caption{The projections of Poincar\'e}
		\label{esfera de poincare}  
	\end{center}
\end{figure}

Let $\,\mathbb{S}^{2} = \{(u,v,w)\in \mathbb{R}^3 \, |\, u^2+v^2+w^2=1\}$ be the 
unit sphere centred at the origin $O$ in $\mathbb R^3$ and identify its tangent 
plane $T_N\mathbb S^2$ at the north pole $N=(0,0,1)$ with the $xy$-plane. Given a 
point ${\bf x} =(x,y,1)\in T_N\mathbb S^2$, the straight line through ${\bf x}$ 
and $O$ intersects $\mathbb S^2$ at the following two points (Fig. \ref{esfera de poincare}):
\begin{equation*}
s_{1}\left( \mathbf{x}\right) =\frac{\mathbf{x}}{\sqrt{1+x^{2}+y^{2}}
}, \,\,\,\, s_{2}\left( \mathbf{x}\right)
=-\frac{\mathbf{x}}{\sqrt{1+x^{2}+y^{2}}}.
\end{equation*}
\medskip

\noindent The maps $s_1:\mathbb R^2\to {\cal H}^+$ and $s_2:\mathbb R^2\to 
{\cal H}^-$ are called {\it the projections of Poincar\'e} where ${\cal H}^+$ (${\cal H}
^-$) denotes the open northern hemisphere of $\mathbb S^2$
$\,\{(u,v,w)\in \mathbb{R}^3 \, |\,\omega > 0\}$ (open southern hemisphere
$\,\{(u,v,w)\in \mathbb{R}^3 \, |\,\omega < 0\}$).	
\vspace{1mm}

\begin{remark}\label{pull-back-form}
The image of each  field  of principal directions $\mathbb{X}_{i}, \, i\in\{1,2\}$ 
under the projection of Poincar\'e $s_j, \, j\in\{1,2\}$ is a direction field diffeomorphic 
to $\mathbb{X}_{i}$ and is the zero loci of the induced quadratic differential form 
$s_{j}^{\ast }\left( \III \right)$.
\end{remark}

\begin{theorem}\label{extendedqde}
The projective extension of the quadratic differential form $\III$ 
is determined by an analytic quadratic differential form $\Phi$ defined on the whole  
unit 2-sphere with the following 
properties:
\begin{itemize}
	\item[i)] the two direction fields defined by $\Phi$ on the upper and lower open hemispheres
	are the zero loci of the induced quadratic forms $s_{1}^{\ast }\left( \III \right)$ 
	and $s_{2} {\ast }\left( \III \right)$,
	\item[ii)] away from the singular points of $\Phi$, the equator is an integral curve of 
	a solution field of $\Phi$ whenever the homogeneous part  $f_n$ of $f$ has no repeated real linear 
	factors.
\end{itemize}
\end{theorem}
\vskip 0.1cm

According to subsection \ref{determ-campos-globalm}, 
the analytic form $\Phi$ referred to in Theorem \ref{extendedqde} defines
globally two direction fields on the unit 2-sphere.

\vskip 0.5cm

\begin{definition} {\rm 
The two direction fields defined by the form $\Phi$ will be denoted $\mathbb{Y}_{1},
\,\mathbb{Y}_{2}$. }
\end{definition}

\noindent {\bf Proof.}  
Consider the map $\varrho : \mathbb{R}
^{3}\backslash \left\{ \omega =0\right\} \rightarrow\mathbb{R}
^{2},\,$ $\left( u,v,\omega \right) \mapsto \left(x,y\right)\,$ where
$\, x = \frac{u}{\omega },\, y = \frac{v}{\omega }.$  The images  
of a pair of antipodal points on the sphere $\mathbb{S}^{2}$ under this map
are the same. We shall now obtain the pullback $\varrho
^{\ast }\left( \III \right) $ of the form $\III$. We rewrite the form $\III$ as

\begin{equation}\label{simplify-princ-eqt} 
\tilde A(x,y) dx^2 + \tilde B(x,y)\ dx dy + \tilde C(x,y) dy^2 =0, 
\end{equation}\vskip 0.2cm 

\noindent where 
$\hskip 1cm\tilde A =  f_{xy} +f_{xy}(f_x)^2 - f_x f_y f_{xx}, \,\, \tilde B =  
f_{yy} \big(1+(f_x)^2\big)-  f_{xx} \big(1+(f_y)^2\big),\,$ \medskip

$\hskip 4.0cm \tilde C =  f_x f_y f_{yy} - f_{xy} - f_{xy} (f_y)^2.$ \medskip

\noindent The replacement of

\hskip 1cm $\begin{pmatrix}
dx & dy%
\end{pmatrix}%
= 
\begin{pmatrix}
\frac{\omega du-ud\omega }{\omega ^{2}} &
\frac{\omega dv-vd\omega }{\omega ^{2}}%
\end{pmatrix}%
= 
\frac{1}{\omega ^{2}}%
\begin{pmatrix}
du & dv & d\omega%
\end{pmatrix}%
\begin{pmatrix}
\omega & 0 \\
0 & \omega \\
-u & -v%
\end{pmatrix}%
$

\noindent in the quadratic form displayed in equation (\ref{simplify-princ-eqt})
leads us to that pullback $\varrho ^{\ast }\left( \III \right) $ is,
\begin{equation}\label{no-simplif-forma}
\frac{1}{\omega ^{4}}%
\begin{pmatrix}
du & dv & d\omega%
\end{pmatrix}%
\begin{pmatrix}
\omega & 0 \\
0 & \omega \\
-u & -v%
\end{pmatrix}%
\begin{pmatrix}
\tilde{A}\left( \frac{u}{\omega },\frac{v}{\omega }\right) & \frac{\tilde{B}}{2}\left( \frac{%
	u}{\omega },\frac{v}{\omega }\right) \\
\frac{\tilde{B}}{2}\left( \frac{u}{\omega },\frac{v}{\omega }\right) & \tilde{C}\left( \frac{%
	u}{\omega },\frac{v}{\omega }\right)%
\end{pmatrix}%
\begin{pmatrix}
\omega & 0 & -u \\
0 & \omega & -v%
\end{pmatrix}%
\begin{pmatrix}
du \\
dv \\
d\omega%
\end{pmatrix},
\end{equation}
where $A, B, C$ are polynomials in $\R[u,v,\om]$ such that
{\small
\begin{eqnarray*}
A(u,v,\om) = \om^{3n-4} \tilde{A}\left( \frac{u}{\omega},\frac{v}{\omega}\right), \,\,\,
B(u,v,\om) = \om^{3n-4} \tilde{B}\left( \frac{u}{\omega },\frac{v}{\omega }\right), \,\,\,
C(u,v,\om) = \om^{3n-4} \tilde{C}\left( \frac{u}{\omega },\frac{v}{\omega }\right),
\end{eqnarray*}}
and
\begin{align}\label{expresiones-ABC}
& A\left( u,v,\omega \right) = \left(F_{uv} (F_u)^2- F_{uu}F_uF_v +\om^{2(n-1)} 
F_{uv}\right)(u,v,\omega), \nonumber \\ 
& B( u,v,\omega) = \left(F_{vv} (F_u)^2- F_{uu}(F_v)^2 +\om^{2(n-1)} (F_{vv}-
F_{uu})\right)(u,v,\omega), \qquad \hskip 2.5cm \nonumber \\ 
& C\left( u,v,\omega \right) = \left(F_{vv}F_uF_v- F_{uv}(F_v)^2 -\om^{2(n-1)} 
F_{uv}\right)(u,v,\omega), \nonumber \\
& F(u,v,\omega) = \sum_{i=0}^{n}\omega ^{n-i} f_i(u,v), \,\,\, 
F_{uu}= \frac{\partial^2 F}{\partial u^2}, \,\,\, F_{uv}= \frac{\partial^2 F}
{\partial u\partial v}, \,\,\, F_{vv}= \frac{\partial^2 F}{\partial v^2}.\qquad
\end{align}
Expanding the product of the three interior matrices of expression (\ref{no-simplif-forma}), 
the pullback $\varrho ^{\ast }\left( \III \right) $ becomes
\begin{equation*}
\frac{1}{\om^{3n}}\
({\begin{matrix}
	du& dv& d\omega%
	\end{matrix}})
\left(\begin{smallmatrix}
\om^2 A & \om^2 \frac{B}{2}& -\om \big(uA+v\frac{B}{2}\big)\\
\om^2 \frac{B}{2}& \om^2 C& -\om \big(u\frac{B}{2}+vC\big)\\
-\om \big(uA+v\frac{B}{2}\big)\ & \ -\om \big(u\frac{B}{2}+vC\big)\ & 
\ u^2A+uvB+v^2C
\end{smallmatrix}\right) 
\left( \begin{matrix}
du \\
dv \\
d\omega%
\end{matrix}\right). 
\end{equation*}
It is worth noting that after multiplying by $\om^{3n}$ and evaluating the last expression 
at $\om =0$ we get the differential form $\, \Big(u^2 A(u,v,0) + uv B(u,v,0) + v^2 C(u,v,0)
\Big) d\om^2$. From this formula it can be seen that the equator is an integral 
solution of both fields away from the singular points if the term
$u^2A+uvB+v^2C$ is not identically zero at $\om =0$. In such a case the form $\III$ 
would be of the type of the positive quadratic forms studied in \cite{Gui}. However, the opposite 
happens since $\om$ is a factor of $u^2A+uvB+v^2C$. To prove it, we will develop 
the polynomial $u^2A(u,v,\om)+uvB(u,v,\om)+v^2C(u,v,\om)$ by doing a recursive application 
of the well-known \medskip

\noindent {\bf Euler's Lemma:}
Let $f\in {\mathbb R}[u,v]$ be a real homogeneous polynomial of degree $n$. Then,
\begin{equation}\label{euler-lemma}
n f(u,v)= u f_u (u,v) + v f_v(u,v).
\end{equation}

\noindent Consider the expressions showed in equation (\ref{expresiones-ABC}). 
For simplicity's sake, in the next steps denote
\begin{equation*}
\Gamma = \om^{2(n-1)}\Big((u^2-v^2) F_{uv} +uv(F_{vv}-F_{uu})\Big). 
\end{equation*}
We thus obtain
\begin{align*}
& u^2 A(u,v,\om) + uv B(u,v,\om) + v^2 C(u,v,\om) = \\
&=- uF_{uu}F_v (uF_u +vF_v) 
+ vF_{vv}F_u(uF_u + vF_v) +F_{uv}(uF_u - vF_v)(uF_u + vF_v) + \Gamma \\
&=\left(uF_u +vF_v\right)
\Big( F_u(vF_{vv}+uF_{uv}) - F_v(uF_{uu}+vF_{uv})
\Big) 
+  \Gamma \quad\\
&=n F
\left( F_u
\bigg( \sum_{i=2}^{n} \om^{n-i}\left[  v\frac{\partial^2 f_i}{\partial v^2}
+ u\frac{\partial^2 f_i}{\partial u\partial v}\right]
\bigg)  - F_v\bigg( \sum_{i=2}^{n} \om^{n-i}\left[ u\frac{\partial^2 f_i}{\partial u^2} 
+ v\frac{\partial^2 f_i}{\partial u\partial v}\right]\bigg)\right)  +  \Gamma.
\end{align*}
By Euler's Lemma, 
$$
v\frac{\partial^2 f_i}{\partial v^2} + u\frac{\partial^2 f_i}{\partial u\partial v} 
= (i-1)\frac{\partial f_i}{\partial v} \qquad \mbox{ and } \qquad
u\frac{\partial^2 f_i}{\partial u^2} + v\frac{\partial^2 f_i}
{\partial u\partial v} = (i-1)\frac{\partial f_i}{\partial u}.
$$
So, after a straightforward calculation we conclude that 
\begin{equation}\label{factorizacion-u2A+uvB+v2C}
 u^2 A(u,v,\om) + uv B(u,v,\om) + v^2 C(u,v,\om) = \om\, T(u,v,\om), 
\end{equation}
where $T \in \R[u,v,\om]$ has the expression 
\begin{align}\label{expresion-T}
T(u, v,\om) &= n f_n \bigg( \frac{\partial f_{n-1}}{\partial u} 
\frac{\partial f_{n}}{\partial v} - \frac{\partial f_{n-1}}{\partial v}
\frac{\partial f_{n}}{\partial u} \bigg)\bigg|_{(u,v)} + \om \Bigg( n f_n (u,v) R(u,v,\omega) 
+ \om^{-2}\Gamma \bigg.
\nonumber\\
& \Bigg. +
\bigg(\sum_{i=1}^{n-1} i\om^{n-i-1} f_i (u,v)\bigg)
\bigg( \om R(u,v,\omega) +
\bigg( \frac{\partial f_{n-1}}{\partial u} 
\frac{\partial f_{n}}{\partial v} - 
\frac{\partial f_{n-1}}{\partial v}\frac{\partial f_{n}}{\partial u}\bigg)\bigg|_{(u,v)} \bigg)\Bigg). 
\end{align}
Similarly $\, R\in \R[u,v,\om]$ is the polynomial 
\begin{align}\label{expresion-R}
&R(u,v,\om) = \bigg(\sum_{i=1}^{n-2} \om^{n-i-2} \frac{\partial f_{i}}{\partial u}\bigg) 
\bigg(\sum_{i=2}^{n-1} (i-1)\om^{n-i} \frac{\partial f_{i}}{\partial v}\bigg)
+ \frac{\partial f_{n}}{\partial u}
\bigg(\sum_{i=2}^{n-2} (i-1)\om^{n-i-2} \frac{\partial f_{i}}{\partial v}\bigg)
\nonumber\\
&+ \frac{\partial f_{n-1}}{\partial u}
\bigg(\sum_{i=2}^{n-1} (i-1)\om^{n-i-1} \frac{\partial f_{i}}{\partial v}\bigg)
- \bigg(\sum_{i=1}^{n-2} \om^{n-i-2} \frac{\partial f_{i}}{\partial v}\bigg) 
\bigg(\sum_{i=2}^{n-1} (i-1)\om^{n-i} \frac{\partial f_{i}}{\partial u}\bigg) 
\nonumber\\
&- \frac{\partial f_{n-1}}{\partial v}\bigg(\sum_{i=2}^{n-1} (i-1)\om^{n-i-1} 
\frac{\partial f_{i}}{\partial u}  \bigg) - \frac{\partial f_{n}}{\partial v}
\bigg(\sum_{i=2}^{n-2} (i-1)\om^{n-i-2} \frac{\partial f_{i}}{\partial u}\bigg)
\nonumber\\
&+ (n-1)\frac{\partial f_{n}}{\partial v}\bigg(\sum_{i=1}^{n-2} \om^{n-i-2} 
\frac{\partial f_{i}}{\partial u}  \bigg) - (n-1)\frac{\partial f_{n}}{\partial u}
\bigg(\sum_{i=1}^{n-2} \om^{n-i-2} \frac{\partial f_{i}}{\partial v}\bigg).
\end{align}
Equality (\ref{factorizacion-u2A+uvB+v2C}) allows us to write the pullback 
$\varrho ^{\ast }\left( \III \right) $ as 

\begin{equation*}\label{preliminar-analityc-form}
\frac{1}{\om^{3n-1}}\
({\begin{matrix}
	du& dv& d\omega%
	\end{matrix}})
\left(\begin{smallmatrix}
\om A & \om \frac{B}{2}& -\big(uA+v\frac{B}{2}\big)\\
\om \frac{B}{2}& \om C& - \big(u\frac{B}{2}+vC\big)\\
- \big(uA+v\frac{B}{2}\big)\ & \ - \big(u\frac{B}{2}+vC\big)\ & 
\ T
\end{smallmatrix}\right) 
\left( \begin{matrix}
du \\
dv \\
d\omega%
\end{matrix}\right). \qquad 
\end{equation*}
After multiplying $\varrho ^{\ast }\left( \III \right) $ by $\,\omega ^{3n-1}\,$ we 
obtain the differential form 

\begin{equation}\label{EDLA}
\Phi \,\, := \,\, \begin{pmatrix}
du& dv& d\omega%
\end{pmatrix}%
\begin{pmatrix}
\omega  A & \omega \frac{B}{2} & 
-\left(u A + v \frac{B}{2}\right) \\
\omega \frac{B}{2} & \omega
C & -\left(u \frac{B}{2} + v C\right) \\
-\left(u A + v \frac{B}{2}\right) & 
-\left(u \frac{B}{2} + v C\right)  & 
T 
\end{pmatrix}
\begin{pmatrix}
du \\
dv \\
d\omega%
\end{pmatrix},
\end{equation}\medskip

\noindent that satisfies the first property. 
To prove the second part of the theorem we need the next

\begin{lemma}\label{express-A+B2}
	The following identities hold
	\begin{align*}
	\left(u A(u,v,\om)+\frac{v B(u,v,\om)}{2}\right)\bigg |_{\om=0} &= 
	\frac{-n v}{2(n-1)} f_n(u,v) \Big|\Hess {f_n}(u,v)\Big|,\\
	\left(\frac{u B(u,v,\om)}{2}+ v C(u,v,\om)\right)\bigg |_{\om=0} &=  
	\frac{n u}{2(n-1)} f_n(u,v) \Big|\Hess {f_n}(u,v)\Big|.
	\end{align*} 
\end{lemma}
When restricting to $\om = 0$ the form $\Phi$ in equation (\ref{EDLA}) we obtain 
by Lemma \ref{express-A+B2} the form
\begin{equation}\label{forma-en-ecuador}
\Phi\big|_{\om =0} := \Big(\frac{n v}{(n-1)} f_n \Big|\Hess {f_n}\Big|\Big) du d\om 
- \Big(\frac{n u}{(n-1)} f_n \Big|\Hess {f_n}\Big|\Big) dv d\om + T(u,v,0) d\om^2.
\end{equation}

\noindent This proves that the equator is locally an integral curve of a direction 
field determined by  $\Phi$.

\begin{remark}\label{Hess-ident-cero}
	If $f_n$ has no repeated real linear factors, then the Hessian polynomial 
	$|\Hess {f_n}|$ is not identically zero. 
\end{remark}
This assertion follows from the classical result ``a binary form $G$ of degree $n$ 
is the $n$th power of a linear form if and only if its Hessian function vanishes 
identically" (Proposition 5.3 of \cite{kung-rota}, Section 3.3.14 of \cite{springer}). \medskip

Remark \ref{Hess-ident-cero} implies that the first two coefficients of the form 
$\Phi\big|_{\om =0}$ of (\ref{forma-en-ecuador}) are polynomials other than the 
zero polynomial provided that $f_n$ 
has no repeated real linear factors. Thus, under this assumption the equator is locally an 
integral curve of only one direction field. This completes the proof of Theorem 
\ref{extendedqde}. \hfill $\Box$ 
\bigskip \medskip	

\noindent {\bf Proof of Lemma \ref{express-A+B2}.}
In the following development it is assumed the notation $\hat f := f_n$. 
Thus,
$$\Big(uA+v\mbox{{\small $\frac{B}{2}$}}\Big)\Big|_{\om=0} = 
\ \frac{1}{2}\Big(- {\hat f}_{uu}{\hat f}_v(u{\hat f}_u +v{\hat f}_v)  
+ ({\hat f}_u)^2(v{\hat f}_{vv}+u{\hat f}_{uv}) + u{\hat f}_{uv}({\hat f}_u)^2 
- u{\hat f}_{uu}{\hat f}_u{\hat f}_v \Big). $$

\noindent By using Euler's formula (\ref{euler-lemma}),
\begin{align*}
2 &\left( uA +v\frac{B}{2}\right)\bigg|_{\om=0} = -n{\hat f}{\hat f}_v{\hat f}_{uu} +(n-1){\hat f}_v
({\hat f}_u)^2 + {\hat f}_{uv}{\hat f}_u(n{\hat f} -v{\hat f}_v) -u{\hat f}_{uu}{\hat f}_u{\hat f}_v  \qquad\\
&= n{\hat f}({\hat f}_u{\hat f}_{uv}-{\hat f}_v{\hat f}_{uu}) + (n-1){\hat f}_v
({\hat f}_u)^2 - {\hat f}_u{\hat f}_v (v {\hat f}_{uv} + u {\hat f}_{uu})\\
&=  n{\hat f}\left( {\hat f}_{uv} \bigg (
\frac{u{\hat f}_{uu}+ v{\hat f}_{uv}}{n-1} \bigg )
- {\hat f}_{uu} \bigg(\frac{u{\hat f}_{uv}+v{\hat f}_{vv}}{n-1}\bigg)
\right) 
+  (n-1){\hat f}_v ({\hat f}_u)^2  - {\hat f}_u{\hat f}_v ((n-1) {\hat f}_u)\\
&=\frac{n}{n-1} {\hat f} \Big( v {\hat f}_{uv}^2 - v {\hat f}_{uu} {\hat f}_{vv} \Big).
\end{align*}
Finally, $\,\displaystyle \left(u A(u,v,\om)+\frac{v B(u,v,\om)}{2}\right)\bigg |_{\om=0} 
= \frac{-n}{2(n-1)} v \Big( f_n(u,v)\Big) \Big|\Hess {f_n}(u,v)\Big|.$	
\vskip 0.3cm

\noindent A similar calculation proves the second identity.
\hfill $\Box$ 

\section{Umbilic Points at Infinity}

The behaviour of the fields $\mathbb{Y}_{1}$ and $\mathbb{Y}_{2}$ restricted 
to the equator of $\mathbb{S}^{2}$ reflects the comportment of the form $\III$ at 
infinity, thus it is relevant to study the singularities of these fields on the equator.

\begin{definition}{\rm 
	A point on the sphere $ \mathbb{S}^{2}$ is a {\it flat point of the 
	differential form $\Phi$} displayed in equation {\rm (\ref{EDLA})} if all the 
coefficients of 
$\Phi$ vanish at that point.}
\end{definition}

The proof of next Lemma follows from equation (\ref{expresion-T}) and Lemma \ref{express-A+B2}. 

\begin{lemma}\label{flat-points-infinity}
	Let $f \in\mathbb{R}\left[x,y\right]$ be a polynomial of degree $n\geq 2$ and
	$p$ be a point on the equator of the sphere $\mathbb{S}^{2}$. Thus $p$ is a 
	flat point of the form $\Phi$ defined in {\rm (\ref{EDLA})} 
	if and only if the homogeneous polynomial $f_n$ vanishes at that point, or 
	both polynomials, $|\Hess f_n|$ and 
	$\nabla f_{n-1}\cdot (\nabla f_n)^\perp$, vanish at $p$. Moreover, when $f_n$ 
	has no repeated real linear factors, $p$ can not be a common zero of $f_n$ and $|\Hess f_n|$.
\end{lemma}

\begin{definition}\label{d:umatinf}{\rm
	A point on the equator of the sphere $\mathbb{S}^{2}$ is an 
	{\em umbilic point at infinity} if it is an isolated flat point of the form $\Phi$.}
\end{definition}

In what follows we need the following notation:
for a polynomial $f \in\mathbb{R}\left[x,y\right] $ of degree $n$ consider its decomposition in homogeneous polynomials, that is,
$\, f(x,y) = \sum_{j=0}^{n} f_{n-j}(x,y)\ $ with $\ f_n(x,y) =  \displaystyle\sum_{\substack{j=0}}^{n-1} a_{j} \, x^{j} y^{n-j},\ 
f_{n-1}(x,y) = \displaystyle\sum_{\substack{ j=0 }}^{n-1} b_{j} \, x^{j} y^{n-1-j} \ $ and 
$\ f_{n-2}(x,y) = \displaystyle\sum_{\substack{ j=0 }}^{n-2} r_{j} \, x^{j} y^{n-2-j}$. From these expressions we define the number 
\begin{equation}\label{number-K}
K := a_{n-2} b_{n-1}^2 + a_{n-1}^2 r_{n-2} - a_{n-1} b_{n-1} b_{n-2}.
\end{equation} 

\begin{remark}
If $p$ is an umbilic point at infinity we can suppose without loss of 
generality that it is the point $(1,0,0)$. Indeed, $p$ can be 
sent to the point  $(1,0,0)$ after a rotation around of the origin of the space with 
coordinates $\{ (u, v, \om)\}$ leaving invariant the $\om$-axis. 	
\end{remark}

The following result provides the Poincar\'e-Hopf index for an umbilic point at infinity and its topological type. We denote by $(u,v)^\perp$ the vector $(v,-u)$ in what follows.

\begin{theorem}\label{indice-umbilic-inf}
Let $f \in\mathbb{R}\left[x,y\right] $ be an $n$-degree polynomial  and suppose that 
the point $p:=(1,0,0)$ is an umbilic point at infinity. If the polynomials $\ \nabla 
f_{n-1}\cdot (\nabla f_n)^\perp$ and $|\Hess f_n|$, have no common real 
linear factors, then
\begin{itemize}
\item[i)] The Poincar\'e-Hopf index of $\mathbb{Y}_k$ at $p$ equals $\frac{1}{2}$,
\item[ii)] $p$ has the topological type of a Lemon if $n=2$; of a Monstar if $n=3$; and 
of a Monstar for $n\geq 4$  whenever $K \neq 0$.	
\item[iii)] $H_f(p) < 0$, where $H_f$ is the homogenization of the 
		Hessian function of $f$ - Definition \ref{d:projhess}.  
	\end{itemize}
\end{theorem}
The proof of Theorem \ref{indice-umbilic-inf} will be given in section \ref{sect-6}. 
\begin{figure}[h!]
	\begin{center}
	\includegraphics[width=7.3cm,height=3.2cm]{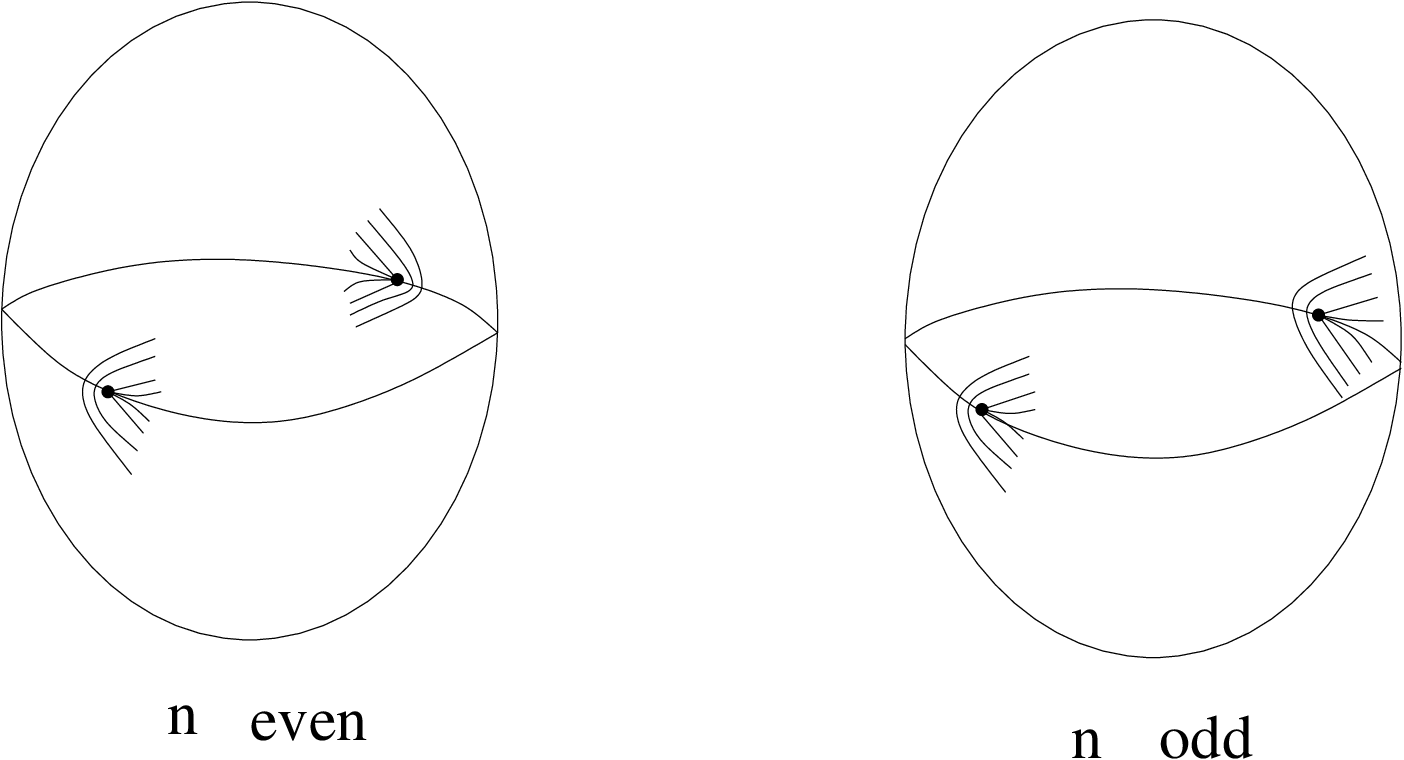}
	\caption{Behaviour of $\mathbb{Y}_k$ at antipodal umbilic points at infinity of Monstar type}  
	\label{antipodales-lemon}
	\end{center}
\end{figure}

In the next result, an upper bound for the number of umbilic points at infinity is 
given.
\begin{theorem}\label{number-umbilic-inf}
	Let $f \in\mathbb{R}\left[x,y\right] $ be a polynomial  of degree $n\geq 2$ 
	and suppose that the $L$ real linear factors of $f_n$ are simple.  If the 
	polynomials $\nabla f_{n-1}\cdot (\nabla f_n)^\perp$ and $|\Hess f_n|$ have 
	$K$ common real linear factors, then the maximal number of umbilic points at infinity is $2L + 2K$. 
\end{theorem}

\noindent{\bf Proof}. 
Note that the number of flat points on the equator of the form $\Phi$ is an upper 
bound for the number of umbilic points at infinity. By Lemma \ref{flat-points-infinity}, 
the number of flat points on the equator is twice the sum of the number of real linear factors of $f_n$
plus the number of common real linear factors of the polynomials $\nabla f_{n-1}\cdot 
(\nabla f_n)^\perp$ and  $|\Hess f_n|$. The second sum is finite because, 
according to Remark \ref{Hess-ident-cero}, the polynomial $\big|\Hess {f_n}\big|$ is not identically zero. 
\hfill $\Box$\medskip

\begin{definition}{\rm
		A homogeneous polynomial on $\R[x, y]$ is called {\it hyperbolic} 
		(resp. {\it{elliptic}}) if its Hessian function has no real linear 
		factors and if it is negative (resp. positive) away from the origin.}
\end{definition}

A better bound is exhibited, for some particular cases, in the next result.

\begin{corollary}\label{ellip-hyperb-cases}
Let $f \in\mathbb{R}\left[x,y\right] $ be a polynomial  of degree $n\geq 2$.  
\begin{itemize}
	\item[i)] If $f_n$ is elliptic, then there are no umbilic points at infinity.
	
	\item[ii)] If $f_n$ is hyperbolic, then the number of umbilic points at 
	infinity is at most twice the number of real linear factors of $f_n$. 	
\end{itemize}
\end{corollary}
\medskip

\noindent{\bf Proof.} 
Note that when $f_n$ is elliptic or hyperbolic, its Hessian polynomial has 
no real linear factors by definition. So, by Theorem \ref{number-umbilic-inf} 
the only polynomial that contributes umbilic points to infinity is $f_n$.
The first assertion follows from the fact that an elliptic homogeneous polynomial  
has no real linear factors (Lemma 3.3 of \cite{GG-OR}). \hfill $\Box$ \medskip

\subsection{Some Remarks about the Homogeneous Case.}\label{remarks-caso-homog}
\begin{remark}\label{ext-homo-case}	
When $f$ is homogeneous, the polynomials $R$ and $T$ displayed in 
{\rm(}\ref{expresion-R}{\rm)} and {\rm(}\ref{expresion-T}{\rm)} respectively, turn out 
	to be $R(u,v,\om) \equiv 0$ and $T(u,v,\om) = \om^{2n-3}\Big((u^2-v^2) f_{uv} 
	+uv(f_{vv}-f_{uu})\Big).$ Indeed, by developing inside the parenthesis of $T$ we have
	$$(u^2-v^2) f_{uv} +uv(f_{vv}-f_{uu}) =  u(u f_{uv}+ v f_{vv}) - v(u f_{uu}+ 
	v f_{uv}) = (n-1) (u f_v - v f_u),$$
	where the last equality is obtained by Euler's\ Lemma. Thus, the analytic 
	differential form $\Phi$ cited in Theorem \ref{extendedqde} becomes 
	\begin{equation}\label{matriz-homo-case}
	\Phi = (\begin{matrix}
	du& dv& d\omega%
	\end{matrix})%
	\left(\begin{smallmatrix}
	\omega  A \qquad\qquad & \omega \frac{B}{2} \quad\qquad & 
	-(u A + v \frac{B}{2}) \qquad\qquad\quad\\
	\omega \frac{B}{2} \qquad\qquad & \omega
	C \qquad\quad & -\left(u \frac{B}{2} + v C\right) \qquad\qquad\\
	-(u A + v \frac{B}{2})\quad & 
	\,\, -\left(u \frac{B}{2} + v C\right)  & 
	\,\,(n-1)\omega^{2n-3}(uf_v - vf_u) 
	\end{smallmatrix}\right)
	\left(\begin{matrix}
	du \\
	dv \\
	d\omega%
	\end{matrix}\right). 
	\end{equation}
\end{remark}

\begin{remark}
The proof of Theorem \ref{indice-umbilic-inf} is valid also in the homogeneous case of degree $n=2$ because it does not depend on the nullity of the coefficients $b_{ij}$.
	For the homogeneous case of degree $n \geq 3$, the origin has the topological type of a Monstar according to Remark 10.1 and Lemma 9.1 of \cite{Gui-2}.
\end{remark}

A better upper bound for the number of umbilic points at infinity is given in the following:

\begin{theorem}\label{number-umb-inf-hom}
	When $f$ is a homogeneous polynomial that has no repeated real linear factors 
	and that is neither elliptic nor hyperbolic, the number of umbilic points at 
	infinity is at most $6n-12$ for $n\geq 5$; 6 if $n = 3$ and $4$ for $n=4$.
\end{theorem}
\medskip 
	
\noindent{\bf Proof.}
If $f$ has exactly $n$ simple real linear factors, then it is a hyperbolic 
polynomial \cite{Guad}, whose situation was analyzed in Corollary \ref{ellip-hyperb-cases}. 
Thus, the maximal number of umbilic points at infinity occurs when $f$ has 
exactly $n-2$ real linear factors and its Hessian polynomial, $2n-4$. In the 
particular case $n=4$, the Hessian polynomial of any homogeneous quartic 
polynomial having exactly two simple distinct real linear factors, has two 
real linear factors whose multiplicity is at least one (see \cite{olver}, 
p. 60).\hfill$\Box$ 
\bigskip

Unlike the general case, it is proven in the next theorem that any flat point on the 
equator is an umbilic point at infinity when $f$ is a homogeneous polynomial.

\begin{theorem}\label{umbilic-inf-hom}
	Let $f$ be a homogeneous polynomial with the property that it has no repeated real 
	linear factors and neither does its  Hessian polynomial. Then, a point is an umbilic point at 
	infinity if and only if it is a flat point on the 
	equator.  	
\end{theorem}

\noindent{\bf Proof.} 
Let $p$ be a point on the equator. After a rotation on the $uv$-plane, we can suppose 
that $p = (1,0,0)$. In a neighborhood of this point the fields $\mathbb{Y}_k$  are 
described by the quadratic differential equation
\begin{equation}\label{carta-u1-homo}
	\begin{pmatrix}
		dv & d\omega%
	\end{pmatrix}%
	\begin{pmatrix}
		\omega\ C(1,v,\omega ) &  -\left(\frac{B}{2} + v C \right)\left(
		1,v,\omega \right) \\ 
		-\left(\frac{B}{2} + v C \right)\left( 1,v,\omega \right) & \quad
		(n-1) \om^{2n-3}\Big( u f_v - v f_u \Big)\left( 1,v\right)%
	\end{pmatrix}%
	\begin{pmatrix}
		dv \\ 
		d\omega%
	\end{pmatrix}%
	=0,  
\end{equation}
whose discriminant, up to a nonzero constant, is

\begin{equation}\label{disci-general}
	\Delta (v,\om)= \Big(B(1,v,\om)+2 v C(1,v,\om)\Big)^2 + 4 (n-1) 
	\om^{2n-2} (vf_u - f_v) C(1,v,\om).
\end{equation}

\noindent Suppose that the origin on the $v\om$-plane is a flat point of 
the differential form given in equation (\ref{carta-u1-homo}). Replacing 
the expressions \medskip

$\,B = f_{vv} (f_u)^2- f_{uu}(f_v)^2 +\om^{2(n-1)} (f_{vv}-f_{uu})\,\,$ and
$\,\,C = f_{vv}f_u f_v- f_{uv}(f_v)^2 -\om^{2(n-1)} f_{uv}\,$ \medskip

\noindent in the discriminant of (\ref{disci-general}) we obtain
\begin{align}\label{discri-adecuado}
	\Delta = &\bigg(\alpha + 2v f_v \Big( f_{vv} f_u - f_{uv} f_v \Big) + 
	\om^{2n-2}\Big(f_{vv} - f_{uu} - 2v f_{uv}\Big)\bigg)^2 \nonumber\\
	& + 4 (n-1) \om^{2n-2}\Big( v f_u - f_v\Big) \Big( \beta - f_{uv} \om^{2n-2}\Big), \qquad 
\end{align}
where  $\alpha$ and $\beta$ are the following polynomials in one variable
$$\alpha(v)= f_{vv} f_u^2 - f_{uu} f_v^2, \quad \beta(v)= f_{vv} f_u f_v - f_{uv} f_v^2. \qquad$$ 

\noindent We remark that point $p$ on the equator is a flat point of 
{\rm (\ref{matriz-homo-case})} if and only if the polynomials, $f\,$ or 
$\,|\Hess {f}|,\,$ vanish at $p$.  Therefore, we consider two cases.
\medskip

\noindent {\bf First case.} Assume that $\,f(1,0) =0$. So, $\,\left( B + 2v C \right)(p) 
= 0$ and $\Delta (0,0) = 0$. It only remains to prove that the origin is an isolated 
singular point. Because $f (1,0) = 0$ and $f$ has no repeated real linear factors, it can be written as
\begin{equation}\label{expresion-f-hom}
	f\left( u,v\right) =v\left( \sum_{\substack{ i=0 }}^{n-1}a_{i} \,
	u^{i} v^{n-1-i}\right),\text{ with }a_{n-1}\neq 0.
\end{equation}
On the one hand, the lowest degree term of $\Delta$ containing only the variable 
$\om$ is $\om^{2n-2}$ whose coefficient is given by the constant part of the single-variable polynomial
$\, 2\alpha (f_{vv} - f_{uu}) - 4(n-1) \beta f_v$. The constant part of 
$\, 2\alpha (f_{vv} - f_{uu})$ is zero. Thus, the coefficient of the monomial 
$\,\om^{2n-2}$ is $\,4 (n-1)^2 a_{n-1}^4\,$ which is positive according to equation (\ref{expresion-f-hom}).
On the other hand, the lowest degree term of $\Delta$ containing only the variable 
$v$ is the monomial $\, n^2 (n-1)^2 a_{n-1}^4 v^2\,$ which appears in the expression 
$\Big(\alpha + 2v f_v ( f_{vv} f_u - f_{uv} f_v )\Big)^2$ and whose coefficient is also positive. 
\medskip
In conclusion, the discriminant $\Delta$ displayed in equation (\ref{discri-adecuado}) 
can be written as 
\begin{equation}\label{discr-hom}
	\Delta (v,\om) = v^2 g_1(v) + \om^{2n-2} g_2(v) + \om^{4n-4} g_3(v),
\end{equation} 
where $g_1, g_2, g_3 \in \R[v]$ are one-variable polynomials such that $g_1(0) > 0$ and $g_2(0) > 0$.
These properties and the parity of the powers appearing in equation (\ref{discr-hom}) 
guarantee that the origin is an isolated singularity.
\bigskip

\noindent {\bf Second case.} Let us suppose that $\,|\Hess f|(1,0) = 0$. 
Since  $f$ has no repeated real linear factor, its Hessian polynomial is not 
identically zero (Remark \ref{Hess-ident-cero}), and $f$ does not vanish at $p$ 
(Lemma \ref{flat-points-infinity}).
Thus, the polynomial $f$ has the expression
\begin{equation}
	f\left( u,v\right) = \sum_{\substack{ i=0 }}^{n} a_{i} \,
	u^{i} v^{n-i},\text{ with }a_{n}\neq 0 \,\, \text{ and } \, a_{n-2} = \frac{(n-1)a_{n-1}^2}{2n a_n}.
\end{equation}
After a straightforward calculation, it is verified that the coefficient 
of the monomial $\,\om^{2n-2}$ is zero and the lowest degree term  in the 
variable $\om$ is the monomial $\,\om^{4n-4}$ whose coefficient is $\,
\left(\frac{n-1}{n a_n}\right)^2 \Big(a_{n-1}^2 - n^2 a_n^2\Big)^2 + 4 
\Big(n-1\Big)^2 a_{n-1}^2 \neq 0$.
The lowest degree term of $\Delta$ including only the variable $v\,$ is the monomial 
$\,\Big(6n^2a_{n-3}a_n^2-(n-1)(n-2)a_{n-1}^3\Big)^2 v^2$. The coefficient of this 
monomial is positive because the polynomial $\,|\Hess f|\,$ has no repeated real linear factors. 
We conclude the proof by noting that the discriminant in this case has the form in equation 
(\ref{discr-hom}) with $g_2 (0) = 0$. \hfill $\Box$ 
\vskip 0.3cm

\section{Umbilic Points on the Finite Part}\label{segunda-sub}

In this section we prove an interesting relation between the indices of all 
the umbillic points of the graph of a real polynomial and the linear factors of the 
highest degree homogeneous part of such a polynomial.

\begin{definition}{\rm 
		We say that a polynomial $f \in \mathbb{R}[x,y]$ is {\it generic} if all of 
		the singular points of the associated fields ${\mathbb Y}_k, 
		k=1,2$, that lie on the equator are umbilic points at infinity.}
\end{definition}

\begin{theorem}\label{poincare-formula}
	Let $f\in \R[x,y]$ be a generic polynomial of degree $n\geq 2$ such that 
	every umbilic point on its graph is isolated. Suppose that $f_n$ has $ L$ real 
	linear factors, and that the polynomials $|\Hess f_n|$ and $\nabla f_{n-1}\cdot 
	(\nabla f_n)^\perp$ have no common real linear factors. Then 
	\begin{equation*}
	\sum_{\substack{$p$ \mbox{ {\rm umbilic}}}} \Ind(p) = 1 - 
	{\textstyle{\frac{1}{2}}}L.
	\end{equation*}
\end{theorem}

The following lemma will be used in the development of the proof.

\begin{lemma}\label{factores-simples}
	Let $f\in \R[x,y]$ be a polynomial of degree $n\geq 2$ and suppose that the polynomials	
	$|\Hess f_n|$ and $\nabla f_{n-1}\cdot (\nabla f_n)^\perp$ have no common real 
	linear factors. Then, every real linear factor of $f_n$ is simple. 
\end{lemma}

\noindent{\bf Proof. }
Suppose that $v$ is a real linear factor of $f_n$. We have thus,

\begin{equation}\label{fn-desarrollada}
f_{n}\left( u,v\right) =v\left( \sum_{\substack{ i=0 }}^{n-1}a_{i} \,
u^{i} v^{n-1-i}\right). 
\end{equation}

\noindent Assume now that $v$ is a repeated linear factor of $f_n$. Thus, $a_{n-1} = 0$, 
which implies that $\,|\Hess f_n| (1,0) = - \big((n-1)\ a_{n-1}\big)^2 = 0$. So, $v$ 
is a factor of the polynomial $\,|\Hess f_n|$. 
The fact that $v^2$ is a factor of $f_n$ leads to $\frac{\partial f_{n}}{\partial  u}\big|_{v=0} = 0$ and $\frac{\partial f_{n}}{\partial  v}\big|_{v=0} = 0$. 
Hence $v$ is also a factor of the polynomial $\nabla f_{n-1}\cdot (\nabla f_n)^\perp$, 
which contradicts the hypothesis. 
Thus, $v$ is a simple real linear factor of $f_n$.\hfill$\Box$

\vskip 0.6cm 

\noindent{\bf Proof of Theorem \ref{poincare-formula}.}
By Lemma \ref{factores-simples}, the $L$ real linear factors of $f_n$ are simple. 
Moreover, each of these factors determines two antipodal points on the equator 
which are flat points of the form $\Phi$. Since every umbilic point on the graph 
of $f$ is isolated by hypothesis, every flat point on the equator is an umbilic point at 
infinity. Therefore, ${\mathbb Y}_k$ has $2L$ umbilic points at infinity.\medskip

Thus, the set of singularities on ${\mathbb S}^2$ of the field ${\mathbb Y}_k$, 
which are all of them isolated, consists of the umbilic points 
at infinity, which lie on the equator ${\mathbb S}^1\subset{\mathbb S}^2$, and the 
finite umbilic points which lie in antipodal pairs on the upper and lower 
hemispheres. \medskip

Applying the Poincar\'e-Hopf Theorem to the direction fields $\mathbb{Y}_k$ defined on 
the whole 2-sphere ${\mathbb S}^2$ which is split into parts to obtain, by Theorem \ref{indice-umbilic-inf},
\begin{align*}
\qquad\qquad\qquad \qquad 2\ &= \sum_{p \in ({\mathbb S}^2\setminus{\mathbb S}^1)} 
\Ind (p)  +  
\sum_{p \in {\mathbb S}^1} \Ind (p)  \\
&= \ 2 \sum_{p\mbox{ umbilic}} \Ind (p) \ + \frac{1}{2} \Big( 2 L \Big). 
\hskip 6.4cm \Box
\end{align*}

\section{Examples}\label{examples}
\begin{figure}[h!] 
\hskip 5cm \includegraphics[width=2.2in]{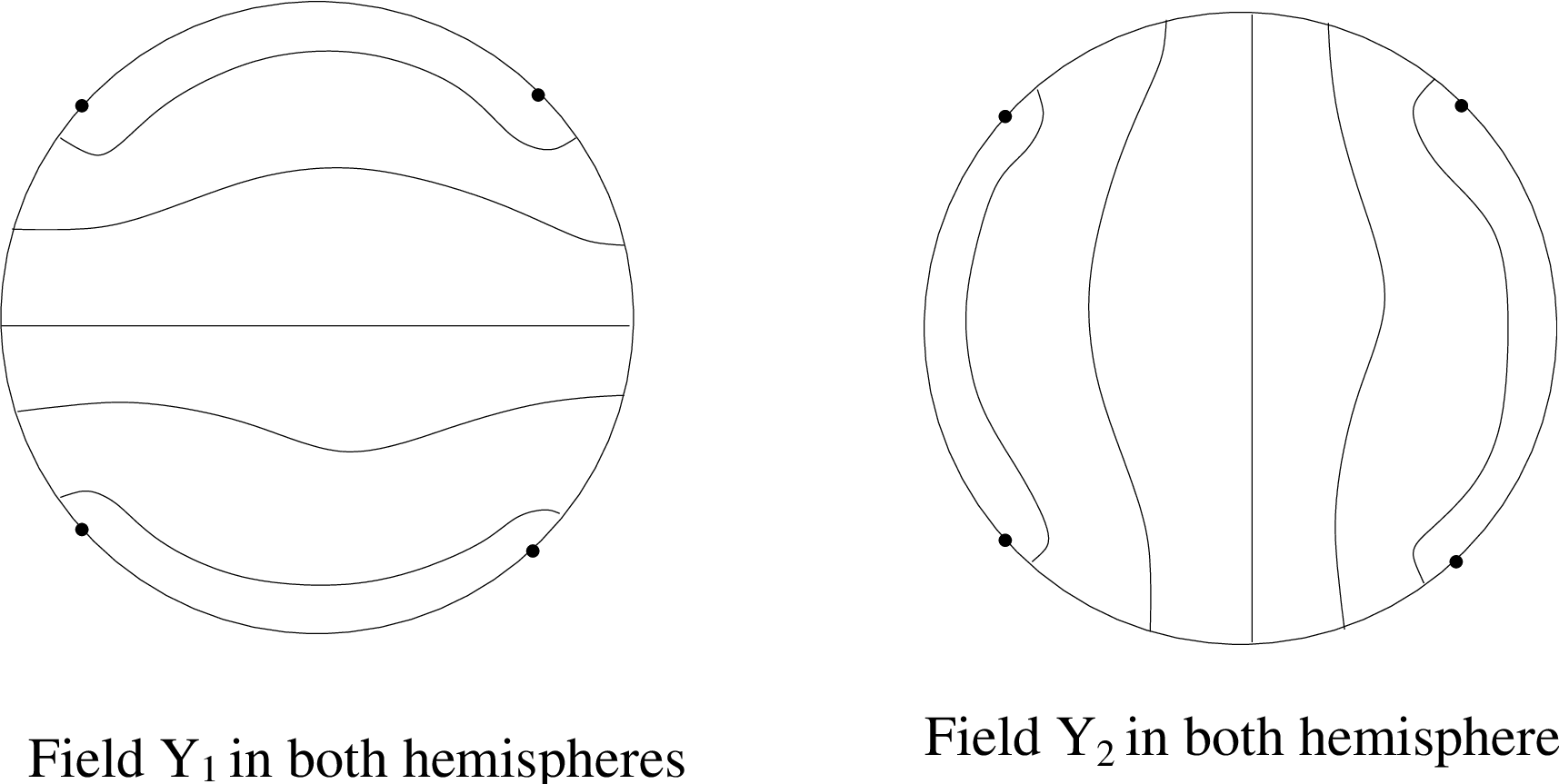}
\caption{Foliations of $\mathbb{Y}_k$ corresponding to  $f = x + 2y + x^2 - y^2$}
\label{silla-no-hom}  
\end{figure}

\begin{examples}\label{ej-cuad-no-hom}
The graph of the polynomial $f(x,y) = x + 2y + x^2 - y^2$ has no umbilic points 
because all of its points are hyperbolic. The homogeneous quadratic part of $f$ is 
a hyperbolic polynomial and each one of its real linear factors gives rise to two  
(antipodal) flat points. Clearly the condition of Theorem \ref{indice-umbilic-inf}, that $ \nabla f_{1}\cdot (\nabla f_2)^\perp$ and $|\Hess f_2|$ have no common real 
linear factors, is satisfied since $|\Hess f_2|= -4$. Since every flat point is isolated, therefore there are four umbilic points at infinity whose topological type is that
of a Lemon. In Fig. \ref{silla-no-hom} we show the foliations of the 
fields $\mathbb{Y}_1, \mathbb{Y}_2$.  
\end{examples}	

\begin{figure}[h!] 
\hskip 5cm \includegraphics[width=2.2in]{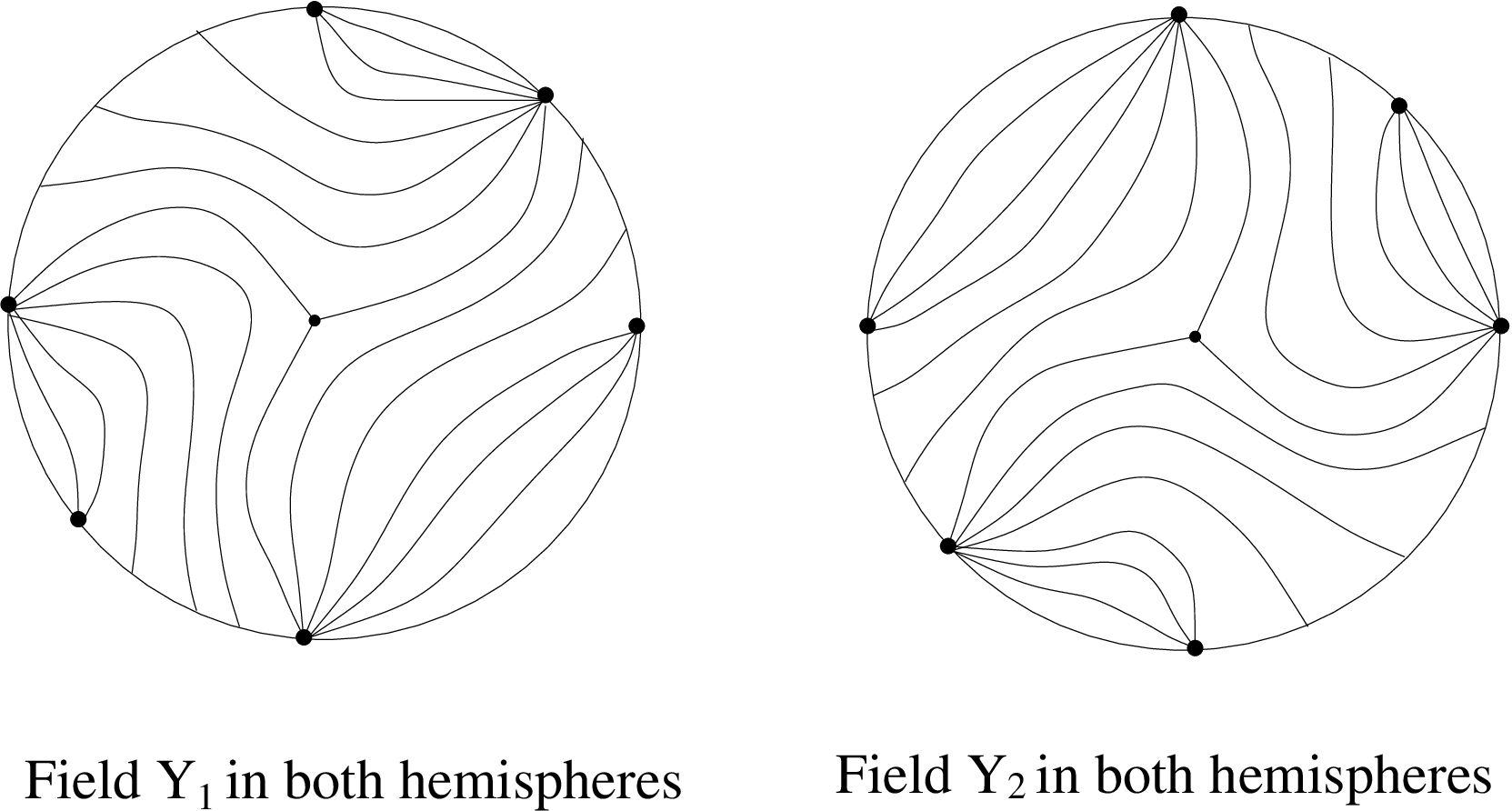}
\caption{Foliations of $\mathbb{Y}_k$ corresponding to $f =  xy + y^2 + xy^2 - x^2y$}
\label{config-grado-tres-no-hom-sfc} 
\end{figure}

\begin{examples}\label{ej-cubico-no-hom}
The graph of the polynomial $f(x,y) =  xy + y^2 + xy^2 - x^2y $ has one umbilic point. 
Since the homogeneous cubic part of $f$ is a hyperbolic polynomial, $|\Hess f_3|$ 
has no real linear factors. Therefore, the condition stated in Theorem 
\ref{indice-umbilic-inf}, that $ \nabla f_{2}\cdot (\nabla 
f_3)^\perp$ and $|\Hess f_3|$ have no common real linear factors, is fulfilled. 
Thus, the flat points on the equator are determined only by the real linear factors 
of $f_3$. Since the umbilic point on the finite part is isolated, every flat point on 
the equator is isolated, which leads to the presence of six umbilic points at
infinity with topological type of a Monstar. In Fig. \ref{config-grado-tres-no-hom-sfc} 
we draw the foliations of the fields $\mathbb{Y}_k$.  
\end{examples}	


\begin{examples}
We provide some examples of homogeneous polynomials that reach the upper bounds 
given in Corollary \ref{ellip-hyperb-cases} and Theorem \ref{number-umb-inf-hom}.
\begin{itemize}
\item For $n\geq 2$, the polynomials $\,f(x,y)= (x^2 + y^2)^n$  are elliptic. In this case there are no umbilic points at infinity because the Hessian polynomial of $f$ is positive away from the origin.
\item
For $n\geq 2$, consider the product of $n$ real linear homogeneous polynomials in 
generic position. They are hyperbolic polynomials of degree $n$ {\rm \cite{Guad}} and 
have $2n$ umbilic points at infinity.\newline
\item In the cubic case, the polynomial $\, f(x,y) = x(x^2 + y^2)$ has 6 umbilic points 
at infinity because the polynomial $\,|\Hess f(x,y)|= 4 (3x^2 - y^2)\,$ has two real linear factors. 
For the quartic case, the polynomial $\, f(x,y) = (x^2 - y^2)(x^2 + y^2)$ reaches the 
bound of 4 umbilic points at infinity. In the remaining cases we do not know if the 
upper bounds are reached.
\end{itemize}
\end{examples}

We now analyze some examples of classical homogeneous quadratic polynomials.

\begin{figure}[ht] 
	\begin{minipage}{0.47\textwidth}
		\qquad\quad\includegraphics[width=2.0in]{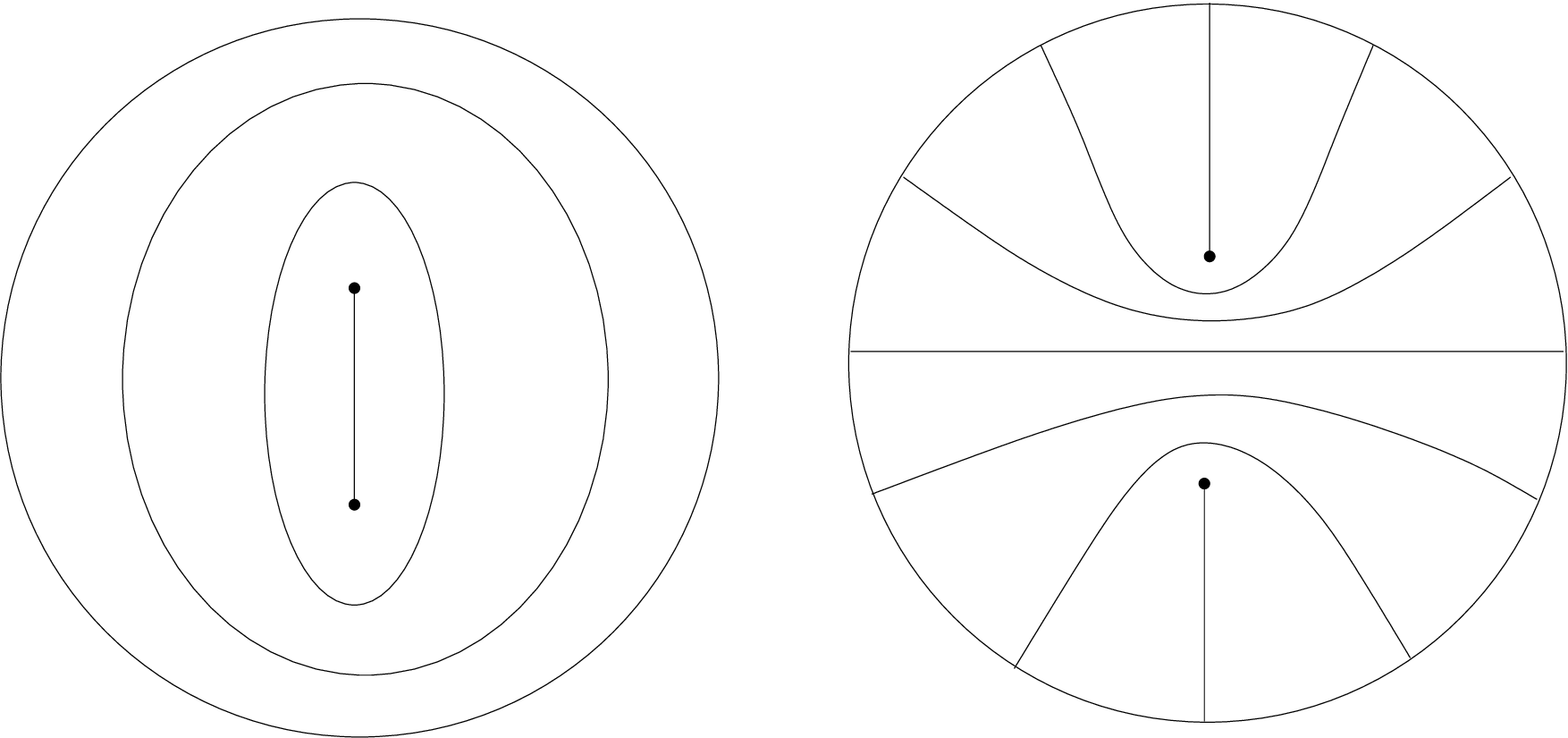}
		\caption{Foliations of $f(x,y) = x^2 + 2 y^2$}
		\label{config-parab-2}  
	\end{minipage}\qquad\,
	\begin{minipage}{0.47\textwidth}
		\qquad\quad\includegraphics[width=2.0in]{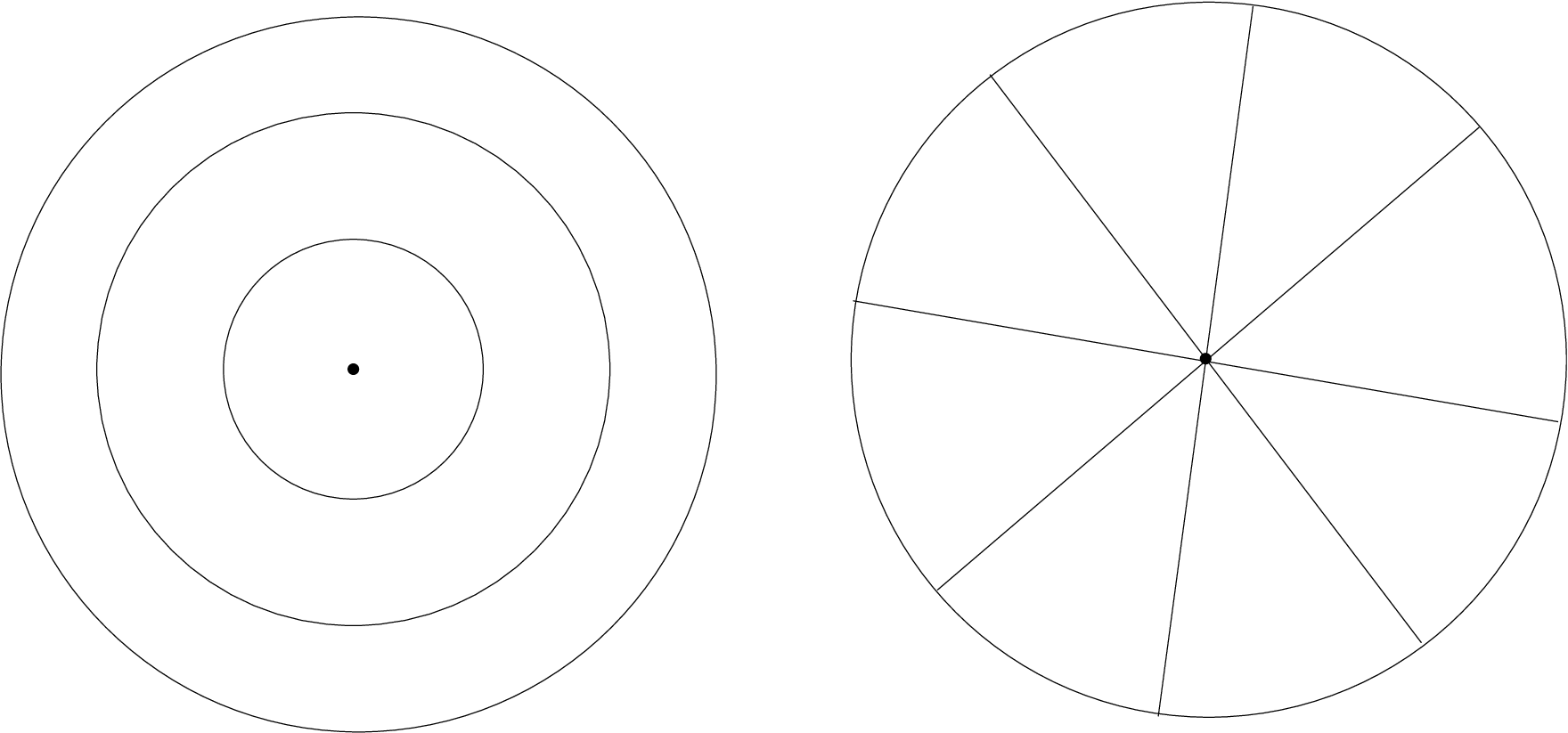}	
		\caption{Foliations of $f(x,y) = x^2 + y^2$}
		\label{config-parab-1}
	\end{minipage}
\end{figure}

\begin{examples} 
{\rm 1)} The graph of $f(x,y) = x^2 + 2 y^2$ has two umbilic points whose 
topological type is a Lemon. Since $f$ is elliptic, by Corollary \ref{ellip-hyperb-cases}, there are no umbilic points 
at infinity, Fig. \ref{config-parab-2}. \vskip 0.3cm
	
\noindent {\rm 2)} The graph of the polynomial $f(x,y) = x^2 + y^2$ has only 
one umbilic point. Because this is a surface of revolution, its meridian and parallel 
curves are lines of principal curvature. There are no umbilic points at infinity 
because $f$ is elliptic,  Corollary \ref{ellip-hyperb-cases}. See Fig. \ref{config-parab-1}.  
\end{examples}

\begin{figure}[h!]
	\begin{center}
		\includegraphics[width=2.05in]{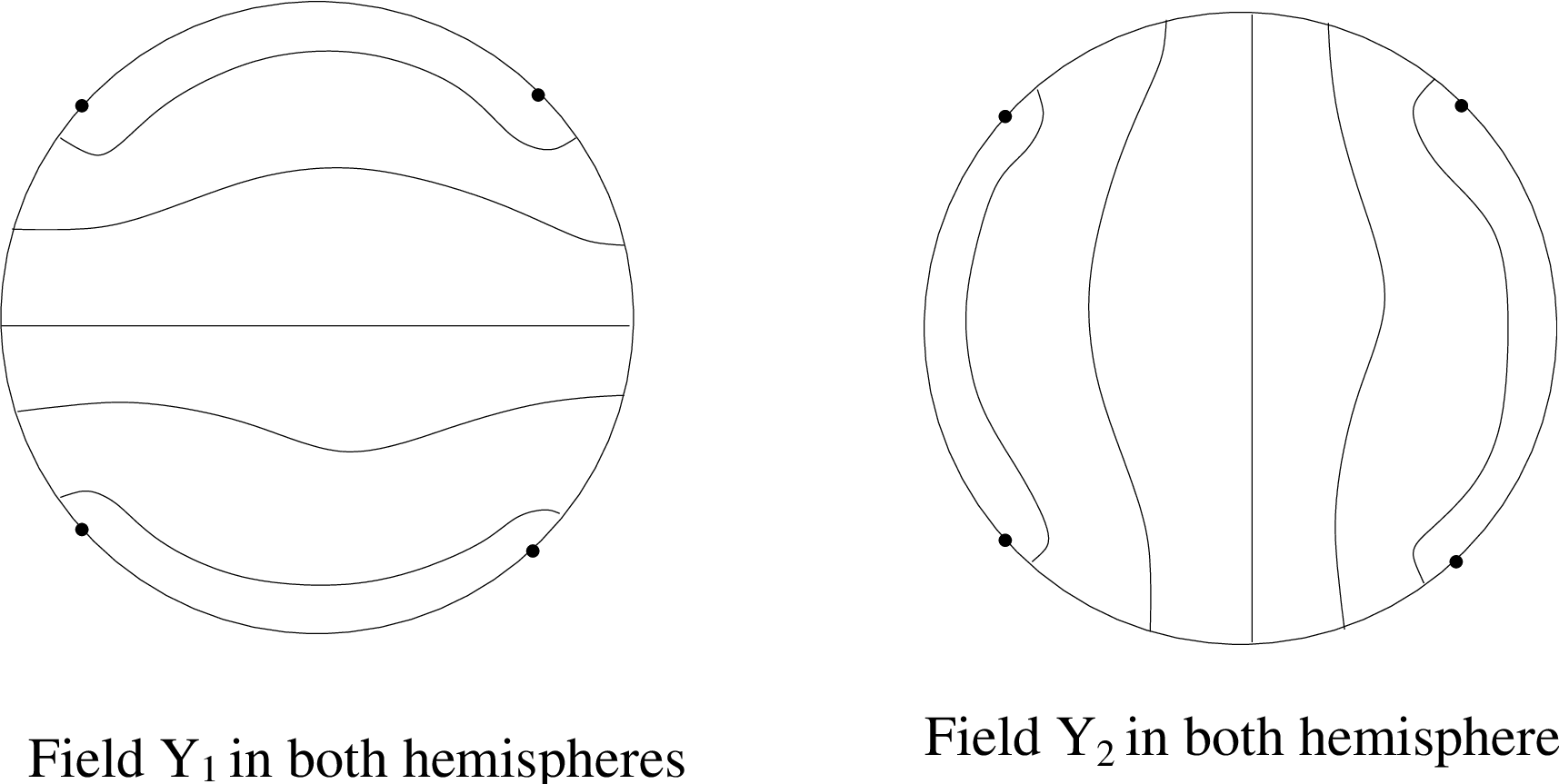}
		\caption{Foliation of $\mathbb{Y}_k$ associated to $f(x,y) = x^2 - y^2$}
		\label{config-silla}
	\end{center}  
\end{figure}

\begin{examples}
The graph of the polynomial $f(x,y) = x^2 - y^2$ has no umbilic points 
because $f$ is a hyperbolic polynomial. Thus, its Hessian polynomial $\,|\Hess f|$ does 
not contribute any flat points on the equator. There are therefore, exactly four flat 
points on the equator determined by $f$, all of which are isolated. In conclusion, 
there are 4 umbilic points at infinity. In Fig. \ref{config-silla}, 
the foliation of the fields $\mathbb{Y}_k, \,k=1,2$ are shown.  
\end{examples}	

\section{Proof of Theorem \ref{indice-umbilic-inf}}\label{sect-6}
Let $q$ be an umbilic point at infinity. According to Lemma \ref{flat-points-infinity}, 
the polynomial $f_n$ vanishes at $q$. So, $v$ is a factor of $f_n$. In accordance with Lemma 
\ref{factores-simples}, $v$ is a simple linear factor of $f_n$. Therefore, $a_{n-1}\neq 0$ 
provided that
$$f_{n}(u,v) = v \left(\sum_{\substack{ i=0 }}^{n-1} a_{i} u^i v^{n-1-i}\right).$$
A simple calculation leads to
$\, H_f(p) = |\Hess f_n| (1,0) = - \big((n-1)\ a_{n-1}\big)^2 < 0.$\medskip
 
Consider now the fields $\mathbb{Y}_{k},\, k=1,2,$ restricted to the set 
$ \{(u,v,\omega)\in \mathbb{S}^2 | u> 0\}$. In the chart $u=1$, they are described by 
the quadratic differential equation

\begin{align}\label{forma-carta-afin}
&\om C(1,v,\omega) dv^2 - \Big(B(1,v,\omega)+2vC(1,v,\omega)\Big) dv d\omega  
+ T(1,v,\omega) d\omega^2 = 0,\qquad\quad \\
\mbox{where}\quad\,\,\, & C(1,v,\omega) = - (n-1) a_{n-1}^3 + \cdots \nonumber\\
& (B+2vC)(1,v,\omega) = - n (n-1) a_{n-1}^3 v - (n-1)(n-2) a_{n-1}^2 b_{n-1} \om + \cdots\nonumber\\
& T(1,v,\omega) = \left\{
\begin{array}{ll}
2 a_{1}^2 b_{1} v + a_1 (1+b_{1}^2) \om + \cdots  & \mbox{ for } n=2 \\
n (n-1) a_{n-1}^2 b_{n-1} v + (n-1)^2 a_{n-1} b_{n-1}^2 \om + \cdots & \mbox{ for } n\geq 3
\end{array} \right. \nonumber \\
\mbox{and \,\,\,} \quad\,\,\, & f_{n-1}(u,v) = \displaystyle\sum_{\substack{ i=0 }}^{n-1} b_{i} \, u^{i} v^{n-1-i}. \nonumber
\end{align}

\noindent We claim that the form displayed in the left-side of equation 
(\ref{forma-carta-afin}) is 
positive in a neighborhood of the origin. Indeed, on the one hand, Theorem 
\ref{extendedqde} guarantees this property in the complement of the equator $\om =0$. On 
the other hand, the restriction to the equator of the discriminant of this form is $\, 
\big( B(1,v,0) + 2v C(1,v,0)\big)^2$, which according to Lemma \ref{express-A+B2} is the 
single-variable polynomial $\Big(\frac{n}{(n-1)} f_n(1,v) |\Hess 
{f_n}|(1,v)\Big)^2$. This polynomial vanishes at a finite number of points since, 
by Remark \ref{Hess-ident-cero}, the polynomial $|\Hess {f_n}(u,v)|$ is different from 
the zero polynomial. 	\medskip

The study  around the origin of the direction fields defined by equation
(\ref{forma-carta-afin}) will be divided into two cases, according to the value of $n$.

\begin{enumerate}[{\bf $i$)}]
\item {\bf Case} $n=2$. In this situation, $a_1 b_1 \neq 0$ and the differential form in 
the left-side of equation (\ref{forma-carta-afin}) becomes
\begin{equation}\label{forma-inf-grado-dos}
\Psi := \Big(- a_{1}^2 \om + \cdots\Big) dv^2 + \Big(2 a_{1}^2 v + \cdots\Big) 
	dv d\omega + \Big(2 a_{1} b_{1} v + (1+b_{1}^2) \om + \cdots\Big) d\omega^2 .
\end{equation} 	
	
To understand the topological behavior of the fields defined by equation
(\ref{forma-inf-grado-dos}) we shall appeal to
the Blowing up method.
\medskip

Consider the correspondence $\gamma:\R^2 \rightarrow \R P^1$ that associates to each point 
on the $v\om$-plane the slope of each straight line defined by the equation 
$$ \Big(- a_{1}^2 \om  + \cdots\Big) \alpha^2 + \Big(2 a_{1}^2 v  + \cdots\Big) 
\alpha \beta + \Big( 2 a_{1} b_{1} v + (1+b_{1}^2) \om + \cdots \Big) \beta^2 =0.$$ 
The image set corresponding to the origin, through $\gamma$, is  $\R P^1$.
The graph of $\gamma$ is a set $\Gamma$ in $\,\R^2\times \R P^1$ which, in 
coordinates, is described as 
\begin{align*}
\Gamma = & \left\{\Big((v,\omega),[\alpha : \beta]\Big) \in \R^2\times \R P^1: 
\big(- a_{1}^2 \om  + \cdots\big) \alpha^2 + \big(2 a_{1}^2 v  + \cdots\big) 
\alpha \beta \right.\\
&\Big. + \big( 2 a_{1} b_{1} v + (1+b_{1}^2) \om + \cdots\big) \beta^2 = 0 \Big\}.
\end{align*}
The discriminant of the form $\Psi$ is $\Delta_\Psi (v,\om) = 4 a_1^4 v^2 + 8 a_1^3 b_1 v \om + 4 a_1^2 (1+b_1^2) \om^2   + \cdots$. Since the Hessian polynomial
of $\Delta_\Psi$ at the origin is $64 a_1^6$, the function $\Delta_\Psi$ has a Morse 
singularity at the origin. Thus (Proposition 2.1, \cite{Bru-Tari}), the set $\Gamma$ is a smooth surface around the 
circle $\{(0,0)\} \times \R P^1$ and the projection $\Pi:M \rightarrow 
\R^2$ defined by $(v,\om,p) \mapsto (v,\om)$, is a local diffeomorphism away 
from the set $\Pi^{-1}(0)$, where $\Pi^{-1}(0) = \{(0,0,p)\}$.  \medskip

Consider the following affine chart on $\R P^1$. Assume $\alpha \neq 0$ and set 
$p=\beta/\alpha$. We define 
$$ F(u,v,p) = \big( 2 a_{1} b_{1} v + (1+b_{1}^2) \om + \cdots \big) p^2 + 
\big(2 a_{1}^2 v + \cdots \big) p + \big( - a_{1}^2 \om + \cdots \big).$$ 
Thus, in the space $\R^3 = \{(v,\om , p)\}$ 
the set $\Gamma$ becomes the smooth surface
$M = \{ (v,\om , p): F(v,\om , p) = 0\}.$

\begin{remark}\label{campo-tg-M}
	The vector field 
	$ \displaystyle \xi =  F_p  \frac{\partial}{\partial v} +  
	 p F_p \frac{\partial}{\partial \om} - 
	\Big(F_v + p F_\om \Big) \frac{\partial}{\partial p}\,$
	is tangent to $M$. Moreover, $\xi$ is a lift on $M$ of the two solution fields 
	{\rm (}\ref{forma-inf-grado-dos}{\rm )}, that is,
	for each point $q\in M$ the vector $\xi(q)$ is sent, under the differential of  
	$\Pi$  into  the vector $F_p \frac{\partial}{\partial v} + p F_p \frac{\partial}{\partial \om}$.
\end{remark}

\begin{proposition}\label{vector-field-lift}
	The vector field $\xi$ has only one zero whose topological type is a saddle.
\end{proposition}

\noindent{\bf Proof.}
The zeros of $\xi$ on $M$ are given by the equations $F = p F_p = F_v + p F_\om =0$. 
The solution set to the system $F =0, \,  p F_p =0$ is the 
set $\{(0,0,p): p\in \R  \}$. Thus, the singular points of $\xi$ are the zeros of 
the cubic single-variable polynomial  $(F_v + p F_\om)\big|_{(0,0,p)} =  p\ Q(p),$ 
where $\,Q(p) = (1+b_{1}^2) p^2 + 2 a_1 b_1 p + a_{1}^2.$ 
Since the discriminant of $Q$ is the negative number 
$-a_{1}^2$, thus the only singular point of $\xi$ is the origin. \medskip

We now prove that $\xi$ has a saddle point at the origin.
Since $\frac{\partial F}{\partial \om}|_{\bar{0}} \neq 0$, the surface
$M$ can be locally written as $\om = g(v,p)$, that is, $F(v,g(v,p),p) \equiv 0$. From this, it follows 
\begin{align}\label{derivadas-impl}
\frac{\partial F}{\partial v} + \frac{\partial F}{\partial \om} \frac{\partial g}{\partial v} = 0, \qquad
\frac{\partial F}{\partial p} + \frac{\partial F}{\partial \om} \frac{\partial g}{\partial p} = 0.
\end{align}

To determine the linear part of $\xi$ at the origin we obtain the linear part 
at the origin of the plane vector field $\bar{\xi} = F_p 
\frac{\partial}{\partial v}  - (F_v + p F_\om) \frac{\partial}{\partial p}$ which
is the projection of $ \xi$ into the $vp$-plane. To accomplish this, write

$$\bar{\xi} = \Big(\alpha_1 v + \alpha_2 p + \cdots\Big) \frac{\partial}{\partial v}  
+ \Big(\beta_1 v + \beta_2 p + \cdots\Big) \frac{\partial}{\partial p}.$$

Note that $\frac{\partial F}{\partial v}|_{\bar{0}} = 0$ because the polynomial 
$F_v + p F_\om$ vanishes at the origin. Using the equalities 
(\ref{derivadas-impl}) we infer that $\frac{\partial g}
{\partial v}\big|_{(0,0)} = 0$; 
and also, $\frac{\partial g}{\partial p}\big|_{(0,0)} = 0$ owing to 
the fact that $F_p(0,0,0) =0$. Thus, \bigskip

\noindent $\alpha_1 = \frac{\partial F_p}{\partial v}\Big|_{(0,0)} = 
\left(\frac{\partial^2 F}{\partial v \partial p} + \frac{\partial^2 
	F}{\partial p \partial \om} \frac{\partial g}{\partial v}\right)\Big|_{(0,0)} 
= \frac{\partial^2 F}{\partial v \partial p}\Big|_{(0,0)} = 2 a_{1}^2 > 0$.

\noindent  $\alpha_2 = \frac{\partial F_p}{\partial p}\Big|_{(0,0)} = 
\left(\frac{\partial^2 F}{\partial p^2} + \frac{\partial^2 
	F}{\partial p \partial \om} \frac{\partial g}{\partial p}\right)\Big|_{(0,0)} 
= \frac{\partial^2 F}{\partial p^2}\Big|_{(0,0)} = 0$.

\noindent  $\beta_1 = -\frac{\partial (F_v + p F_\om)}{\partial v}\Big|_{(0,0)} 
= - \frac{\partial^2 F}{\partial v^2}\Big|_{(0,0)} = 0$. 

\noindent  $\beta_2 = -\frac{\partial (F_v + p F_\om)}{\partial p}\Big|_{(0,0)} = -  
a_{1}^2 < 0$. \medskip

\noindent We conclude that the origin is a saddle point (Fig. \ref{surface-espacio-jets}).
\begin{remark}\label{pto-silla-proy-lemon}{\rm (\cite{Bru-Fid}, p.152)}
	The projection into the $v\om$-plane of the integral curves of the field $\xi$ 
	under the map $\Pi$ leads to the conclusion that the origin is a singular point 
	whose topological type is a Lemon {\rm (Fig. \ref{lemon-point})}. 
\end{remark}
 
\begin{figure}[ht] 
	\begin{minipage}{0.47\textwidth}
		\qquad\includegraphics[width=1.8in]{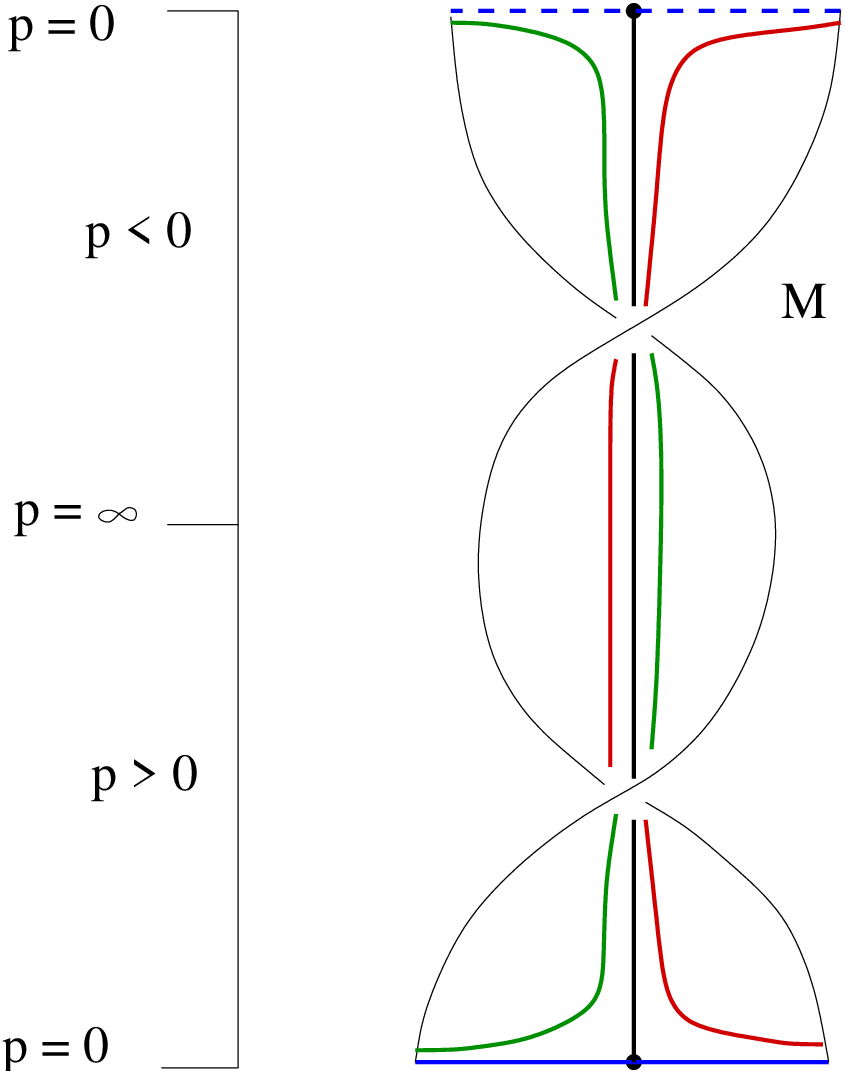}
		\caption{One saddle point in the surface $M$}
		\label{surface-espacio-jets} 
	\end{minipage}\quad
	\begin{minipage}{0.44\textwidth}
		\qquad\includegraphics[width=2.3in]{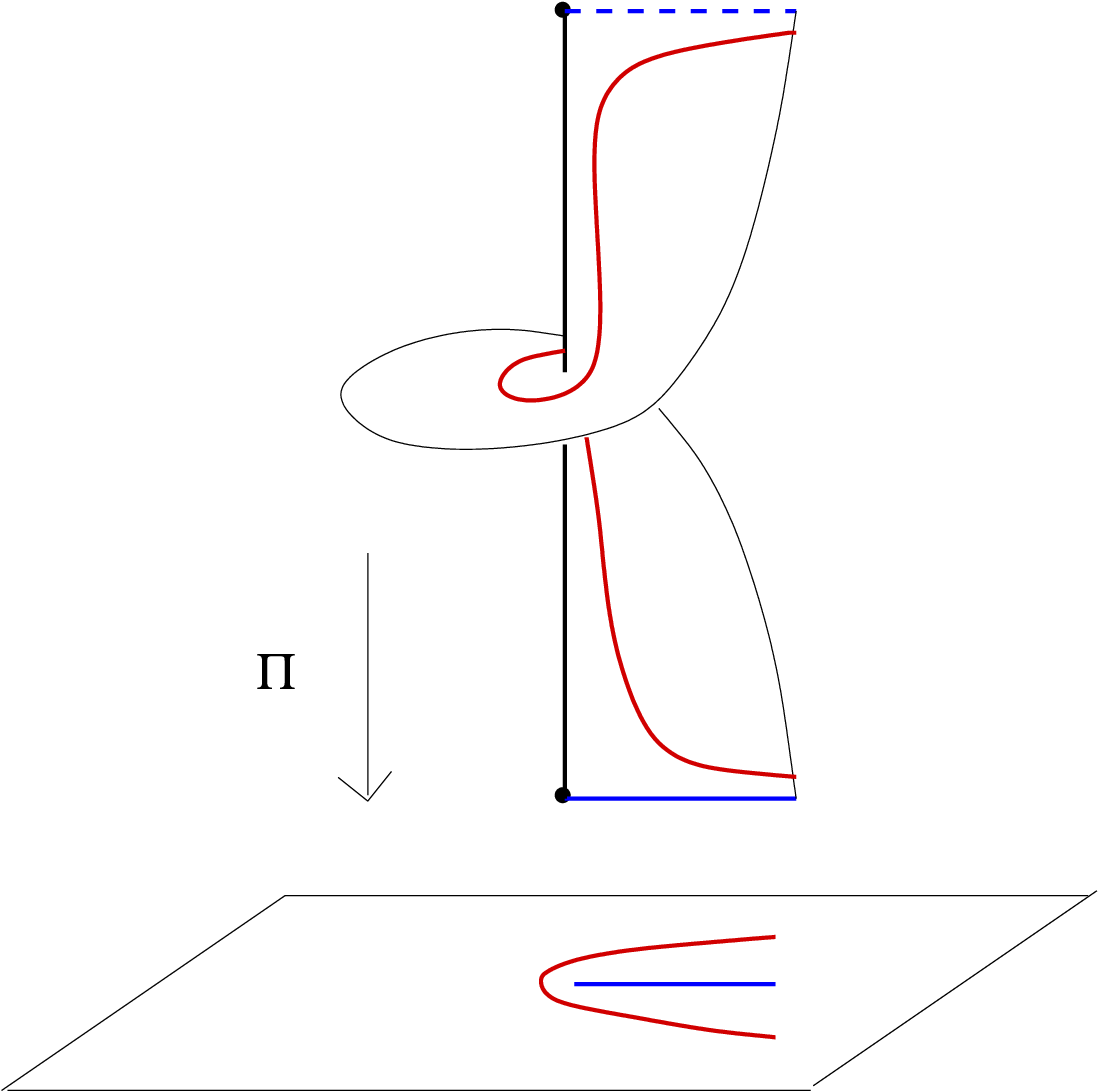}	
		\caption{Projection of a saddle point}
		\label{lemon-point}
	\end{minipage}
	\qquad\qquad
\end{figure}
\item {\bf Case} $n \geq 3$. The differential form (\ref{forma-carta-afin}) 
becomes, after dividing it by $(n-1)a_{n-1}$, into 
\begin{align}\label{forma-inf-grado-general}
\Psi := &\Big(- a_{n-1}^2 \om + \cdots\Big) dv^2 + \Big(n a_{n-1}^2 v 
+ (n-2) a_{n-1} b_{n-1} \om + \cdots\Big) dv d\omega \qquad \nonumber\\ 
&+ \Big(n a_{n-1} b_{n-1} v + (n-1) b_{n-1}^2 \om 
+ \cdots\Big) d\omega^2 .
\end{align}
Since the quadratic part of the discriminant of $\Psi$ is $ n^2 a_{n-1}^2 
(a_{n-1} v + b_{n-1} w)^2$, the origin is not a Morse singularity of $\Delta_\Psi$. 
Thus, we can not proceed as in the previous case. \medskip

As a first step, consider the change of coordinates on the $v\om$-plane
$$ \begin{pmatrix}
v \\
\om 
\end{pmatrix} = 
\begin{pmatrix} 
1 &  - \frac{b_{n-1}}{a_{n-1}}\\
0 &  1
\end{pmatrix} 
\begin{pmatrix} 
X \\ Y \end{pmatrix}.$$
The differential form $\Psi$ of (\ref{forma-inf-grado-general}) is transformed into the differential 
form
\begin{align*}
\widetilde{\Omega} =\widetilde{A}(X,Y) dX^2 + \widetilde{B}(X,Y) dX dY + \widetilde{C}(X,Y) dY^2,
\end{align*}
where 
\begin{align*}
\widetilde{A}(X,Y) &= a\left(X - \frac{b_{n-1}}{a_{n-1}} Y,Y\right) = Y (- a_{n-1}^2 + \cdots),\\
\widetilde{B}(X,Y) &= \Big( b-2 \frac{b_{n-1}}{a_{n-1}} a
\Big)\left(X - \frac{b_{n-1}}{a_{n-1}} Y, Y\right) = n a_{n-1}^2 X + \cdots,\\
\widetilde{C}(X,Y) &=
\left( \frac{b_{n-1}^2}{a_{n-1}^2} a -\frac{b_{n-1}}{a_{n-1}} b + c
\right)\left(X - \frac{b_{n-1}}{a_{n-1}} Y, Y\right) =  \cdots,
\end{align*}
and $a, b, c$ are the first, second and third coefficients of 
(\ref{forma-inf-grado-general}). In what follows we consider the differential
form $\Omega$ obtained of dividing $\widetilde{\Omega}$ by $a_{n-1}^2$, that is,
\begin{align}\label{forma-nueva}
\Omega := A(X,Y) dX^2 + B(X,Y) dX dY + C(X,Y) dY^2,
\end{align}
$A(X,Y) = -Y + \displaystyle\sum_{i+j=2}^{3n-3} a_{ij} X^i Y^j, \ B(X,Y) = nX + \sum_{i+j=2}^{3n-3} b_{ij} X^i Y^j, \ C(X,Y) = \sum_{i+j=2}^{3n-3} c_{ij} X^i Y^j$.\medskip

As some terms of the previous coefficients depend on the degree $n$ we split the following analysis into two cases. \medskip 

\begin{enumerate}[{\bf a)}]
	\item {\bf Case} $n=3$. Consider the function $G_{\Omega}:\R^2 
\rightarrow \R^3$ that associates to each pair $(X,Y) \in \R^2$ the coefficients 
$\, (C(X,Y),\, B(X,Y), \, A(X,Y))\, $ of the differential 
form $\Omega$. The Jacobian matrix D$G_\Omega $ of the map $G_\Omega$ at the origin is 
$\begin{pmatrix} 
0 & 0\\
3 & 0\\
0 & -1
\end{pmatrix}. $ 

Since the rank of D$G_\Omega$ at the origin is 2, the origin is called a {\it semi-simple singular point of $\Omega$} according to the notation of \cite{Gui-2}. Moreover, $\Omega$ has the type of $E(\lambda)$ for $\lambda =3$. By Remark 10.1 and Lemma 9.1 of \cite{Gui-2}, the origin is a singular point of $\Omega$ with topological type of a Monstar. \bigskip 

\item Case $n\geq 4$.  As a second step, we will now transform the form (\ref{forma-nueva}) into a 
suitable differential form through the Blowing up method.
On the $VW$-plane consider the isomorphism 
$$\phi:\R^2\setminus\{V=0\} 
\rightarrow \R^2\setminus\{X = 0\} \mbox{ defined as } \phi (V,W) = (V, VW) = (X,Y),$$ 
that is, $V = X, W= Y/X$.	The pullback $\phi^* \Omega$ of $\Omega$ is the 
differential form
\begin{equation}
\phi^* \Omega = A dV^2  + B dV dW + C dW^2
\end{equation}
where $  a,b,c$ denote respectively the first, second and third coefficients of $\Omega$, and
\begin{align*}
A (V,W)&= W^2 c(V, VW) + W b(V, VW) + a(V, VW), \\
B (V,W)&= V b(V, VW) + 2 V W c(V, VW),\\
C (V,W)&= V^2 c(V, VW). 
\end{align*}
Since $V$ is a factor of $a, b, c$, rewrite the form $\phi^* \Omega$ as $ \phi^* \Omega = V \Omega_1$ where 
$$ \Omega_1 := A_1 \ dV^2 + V B_1 \ dV dW +  V^2 C_1 \ dW^2,$$
and $A_1, B_1, C_1$ are polynomials in $\R[V,W]$.
A straightforward calculation shows that 
\begin{align}\label{expres-A1}
A_1 (V,W) =  (n-1) W + V W g_1,\quad 
B_1 (V,W) = n + V g_2, \quad
C_1 (V,W) =  V g_3,
\end{align}
where $g_i, \, i\in \{1,2,3\}$ is a polynomial in $\R[V,W]$. Note that the origin is the only singular point  of the form $\Omega_1$ on the line $V=0$. \medskip 

In a neighborhood of the origin on the $VW$-plane the two fields of directions defined by $\Omega_1$ are described by 
$$ \Big(-2 A_1\Big) \, dV + \Big(- V B_1 + (-1)^{k} \sqrt{V^2 (B_1^2 - 4 A_1 C_1)} \Big) dW = 0, \quad k\in \{1,2\}.$$

Let's denote by ${\cal F}_k(\Omega_1), k \in \{1,2\}$ the foliations corresponding to these direction fields. 
Consider now the vector fields 
$$ Y_k(V, W) := (V\, T_k , 2 A_1), \quad \mbox{with}\quad
T_k = - B_1 + (-1)^k \sqrt{B_1^2 - 4 A_1 C_1}.$$

In a punctured neighborhood of the origin the foliation ${\cal F}_1(\Omega_1)$ is 
tangent to the vector field $Y_1$ if $V > 0$, and tangent to the vector field $Y_2$ 
if $V < 0$. Analogously, the foliation ${\cal F}_2(\Omega_1)$ is tangent to $Y_2$ when 
$V > 0$, and tangent to $Y_1$ for $V < 0$.\medskip

From expressions (\ref{expres-A1}) we infer that $T_1(0,0) = - 2 n,\ $  $T_2(0,0) = 0$ and the linear part of $Y_1$ at the origin is
\begin{eqnarray*}
	DY_{1}|_{(0,0)} &=&\left.\left( 
	\begin{array}{cc}
		V \frac{\partial }{\partial V }T_{1} + T_1  & 
		V \frac{\partial }{\partial W }T_{1} \\	
		2\frac{\partial }{\partial V} A_1 \qquad& 2\frac{\partial }{%
			\partial W } A_1 
	\end{array}%
	\right) \right\vert _{\left( 0,0\right) } 
	= \left( 
	\begin{array}{cc}
	-2n   & 0 \\	
	0 \,\, & n-1 
\end{array}%
\right).
\end{eqnarray*}
Therefore, the origin is a saddle point of the field $Y_1$ and the eigenspaces of $DY_{1}|_{(0,0)}$ are the coordinate axes.\medskip

On the other hand, consider the vector field 
$$\ Z_1 (V,W) = (2 V C_1(V,W),\ T_1 (V,W)).$$
Since $T_{1}(0,0)\neq 0$, the origin is a nonsingular point of this field. Moreover, 
because of the equality $T_1(V,W)\, T_2 (V,W) = 4 A_1 (V,W)\, C_1 (V,W)$, the field $Z_1$ satisfies  the relation
$$ 2 A_1 (V,W)\, Z_1 (V,W) = T_1 (V,W) \, Y_2 (V,W),$$
which proves that $Z_1$ is tangent to the foliation of $Y_2$ (Fig. \ref{fig-12}).

\begin{figure}[h!]
	\begin{center}
		\includegraphics[width=3.6cm,height=2.3cm]{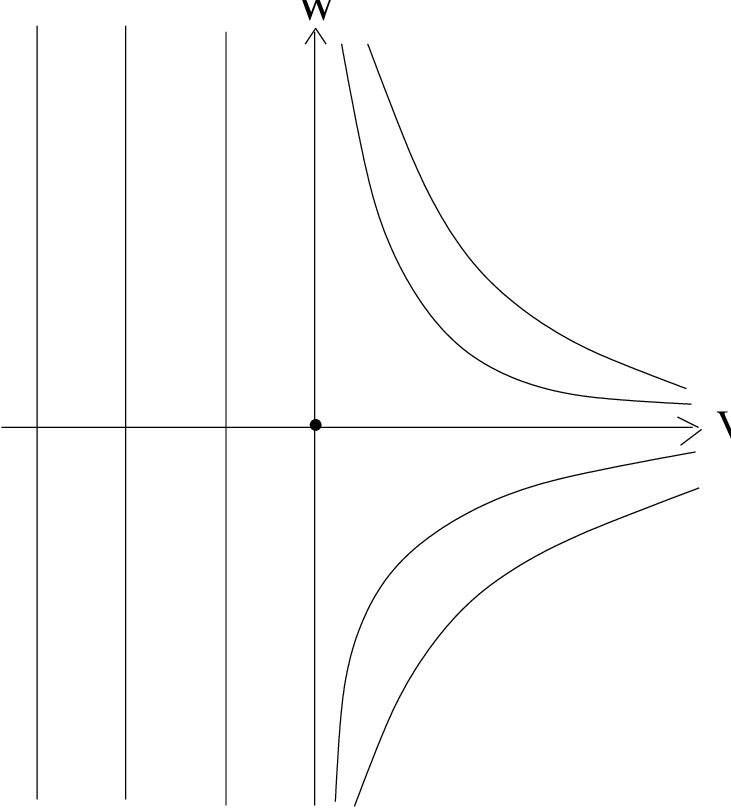}
		\caption{Foliation ${\cal F}_1(\Omega_1)$}
		\label{fig-12}
	\end{center}  
\end{figure}

\bigskip 

In order to carry out a complete analysis of the singularity we will do another 
bluwing up of $\Omega$ as a third step.

Consider the map $\varphi (V,W) = (VW^2, W) = (X,Y)$. The pullback $\varphi^*
\Omega$ of the form $\Omega$ is
$$\varphi^*\Omega = (W^4 A)\ dV^2 + W^2 (B + 4VWA)\ dV dW + (4V^2W^2 A + 2VW B + C)\ dW^2.$$

A straightforward calculation shows that $\varphi^*\Omega = W^3\, \Omega_2$ where
\begin{align}\label{omega2}
\Omega_2 = W^2 A_2 \ dV^2 + W B_2 \ dV dW + C_2 \ dW^2,
\end{align}
where
$\quad  A_2 (V,W) = -1 + W h_1,$ \quad $B_2 (V,W) = b_{02} + (n-4) V + W h_2,$ \medskip 

\hskip 1.5cm $C_2 (V,W) = c_{03} + (c_{11} + 2 b_{02}) V + (2n-4) V^2 + W h_3,\, \ $ and \medskip

\hskip 1.5cm$ h_k,\ $ is a polynomial in $\R[V,W]\ $ for $\ k \in \{1,2,3\}$.
\medskip

The two fields of directions defined by $\Omega_2$ are described, in a neighborhood of the origin on the $VW$-plane by
$$ \Big(- W B_2 + (-1)^{k} \sqrt{W^2 (B_2^2 - 4 A_2 C_2)}\Big) \, dV + \Big(-2 C_2 \Big) dW = 0, \quad k\in \{1,2\}.$$

Let's denote by
${\cal F}_k(\Omega_2), \, k \in \{1,2\}$ the foliations of these fields. Consider the vector fields 
$$Y_k(V, W) := (2 C_2,\, W\, T_k), \quad \mbox{with}\quad
T_k = - B_2 + (-1)^k \sqrt{B_2^2 - 4 A_2 C_2}.$$

In a punctured neighborhood of the origin the foliation ${\cal F}_1(\Omega_2)$ 
is tangent to the vector field $Y_1\,$ if $\, W > 0$, and tangent to the vector 
field $Y_2\,$ if $\, W < 0$. Analogously, the foliation ${\cal F}_2(\Omega_2)$ is 
tangent to $Y_2$ when $W > 0$, and tangent to $Y_1$ for $W < 0$.\medskip

On the other hand, in the following analysis we will also consider the vector fields 
$$\ Z_k (V,W) = (2\, V\, C_2(V,W),\ T_k (V,W)) \  \mbox{ for }  k\in \{1,2\}.$$
Because of the equality $T_1(V,W)\, T_2 (V,W) = 4 A_2 (V,W)\, C_2 (V,W)$, the field $Z_k,\ k\in\{1,2\}$ satisfies  the relation
\begin{equation}\label{relat-zk-yk-grado3}
2\, C_2 (V,W)\, Z_k (V,W) = T_k (V,W) \, Y_{3-k} (V,W).
\end{equation}

In what follows we will need the following coefficients:  
\begin{align}\label{expresions-coef}
b_{02} = \frac{(n-2) (n-3) K}{(n-1) a_{n-1}^3}, \quad 
c_{11} = \frac{2 n (n-2) K}{(n-1) a_{n-1}^3}, \quad  
c_{03} = \frac{ 2 (n-2)^2 K^2}{(n-1) a_{n-1}^6}, \quad 
\end{align}
where $K := a_{n-2} b_{n-1}^2 + a_{n-1}^2 r_{n-2} - a_{n-1} b_{n-1} b_{n-2} \ $
and $\ f_{n-2} (u,v) = \displaystyle\sum_{\substack{ i=0 }}^{n-2} r_{i} \, u^{i} v^{n-1-i}.$

The singular points of $\Omega_2$ on the line $W=0$ are given by the equation $\ C_2 
(V,0) = 0,$ that is, $\ c_{03} + (c_{11} + 2 b_{02}) V + (2n-4) V^2$. The
discriminant of this cuadratic equation is $\Delta = \frac{4 (n-2)^2 K^2}{(n-1)^2 a_{n-1}^6}$ which is positive if and only if $K \neq 0$. 
Therefore, the singular points on $W=0$ are
$$p_j := \left(-\frac{(2n-3) K\ }{2 (n-1) a_{n-1}^3} (-1)^j \sqrt{\frac{K^2}{4(n-1)^2 a_{n-1}^6}}\,, \ 0\right), \ \mbox{ for } j \in \{1,2\}.$$

Suppose $K/a_{n-1} >0$. This condition implies that $\, V_1 = -\frac{K}{a_{n-1}^3},\
V_2 = -\frac{(n-2) K}{(n-1) a_{n-1}^3}\, $ and both roots are negative.
From expressions (\ref{expresions-coef}) we derive that $\ T_1(p_1) = -
\frac{4 K }{(n-1) a_{n-1}^3}, $ $\, T_1(p_2) = -\frac{2 (n-2) K }{(n-1) a_{n-1}^3}\ $  and the linear part of $Y_1$ at $p_j$ is
\begin{eqnarray*}
DY_{1}|_{P_j}  = \left\{
\begin{aligned}
&\left( 
\begin{array}{cc}
	-\frac{4(n-2) K}{(n-1) a_{n-1}^3} & 2\frac{\partial  C_2}{\partial W }(p_1) \\	
	0 \,\, & -\frac{4 K}{(n-1) a_{n-1}^3}
\end{array}%
\right)	\ \mbox{ for } j=1, \\
&\left( 
\begin{array}{cc}
\frac{4(n-2) K}{(n-1) a_{n-1}^3}   & \ \ 2\frac{\partial  C_2}{\partial W }(p_1)\\	
	0 \,\, & \ \ -\frac{2 (n-2) K }{(n-1) a_{n-1}^3}
\end{array}%
\right) \ \mbox{ for } j=2.
\end{aligned}\right.
\end{eqnarray*}

Therefore, $p_1$ is a node point of the field $Y_1$, and $p_2$ is a saddle point. On 
the other hand, since $T_{1}(p_j)\neq 0 \ $ for $j\in\{1,2\}$, the point $p_j$ is a 
nonsingular point of $Z_1$. Moreover, according to (\ref{relat-zk-yk-grado3}) the field $Z_1$ is tangent to the foliation of $Y_2$. We remark that analogous results are
obtained when $K/a_{n-1} <0$. \medskip 

From all the previous analysis we conclude that the the origin on the $v\om$-plane has the topological type of a Monstar when $K \neq 0$.
This completes the proof of Theorem \ref{indice-umbilic-inf}. 
\hfill $\Box$
\end{enumerate} 
\end{enumerate}

\vskip 0.2cm
\noindent {\bf Acknowledgments}
The second author thanks to the PASPA-DGAPA, UNAM Program 
for its finantial support during the preparation of this work and also 
thanks to the Institute of Technology, Tralee for its hospitality.

\end{document}